\documentclass{amsart}
\usepackage{amssymb,amsmath}
\newcommand{\RR}{{\mathbb R}}

\newcommand{\e}{\varepsilon}

\newcommand{\Lap}{\Delta}
\newcommand{\de}{\delta}
\newcommand{\del}{\partial}
\newcommand{\om}{\omega}
\newcommand{\Om}{\Omega}

\newcommand{\si}{\sigma}

\newcommand{{\loc}}{{\ell\mathrm oc}}

\def\meanint{{\diagup\hskip -.42cm\int}}

\newtheorem{theorem}{Theorem}
\newtheorem{lemma}{Lemma}
\newtheorem{proposition}{Proposition}
\newtheorem{corollary}{Corollary}

\newtheorem{remark}{Remark}

\begin{document}

\title[Differentiability of Solutions]
{Differentiability of Solutions to Second-Order Elliptic Equations via Dynamical Systems}
\author{Vladimir Maz'ya}
\address{Link\"oping University,  University of Liverpool}
\thanks{The first author was partially supported by the UK Engineering and Physical Sciences Research
Council  grant EP/F005563/1.}

\author{ Robert McOwen}
\address{Northeastern University}
\date{July 12, 2010}
\maketitle

\begin{abstract}
For a second-order elliptic equation in divergence form we investigate conditions on the coefficients which imply that all solutions are  Lipschitz continuous or differentiable at a given point. We assume the coefficients have modulus of continuity satisfying the square-Dini condition, and obtain additional conditions that examples show are sharp. Our results extend those of previous authors who assume the modulus of continuity satisfies the Dini condition. 
 Our method involves the study  of asymptotic properties of solutions to a dynamical system that is derived from the coefficients of the elliptic equation. 


\end{abstract}

\addtocounter{section}{-1}
\section{Introduction}
\smallskip\noindent

We consider the regularity of  weak solutions of a linear uniformly elliptic equation in divergence form
in an open set $U$ of $\RR^n$ for $n\geq 2$:
\begin{equation}
{\mathcal L}u:=\del_i(a_{ij}(x)\del_j u)=0 \quad\hbox{in}\  U,
\label{eq:Lu=0}
\end{equation}
where we have used the summation convention, the $a_{ij}=a_{ji}$ are bounded, measurable, real-valued functions, and
by a {\it weak solution} of (1) we mean that $u\in H_{\loc}^{1,2}(U)$, i.e.\ $\nabla u$ is locally square-integrable,  and satisfies
\begin{equation}
\int_U a_{ij}(x)\,\del_j u\,\del_i\eta\,dx=0\quad\hbox{for all $\eta\in C_0^\infty(U)$.}
\label{def:Lu=0}
\end{equation}
The classical results of De Giorgi \cite{DG} and Nash \cite{N} show that $u$ is locally H\"older continuous in $U$. When the coefficients are continuous in $U$, then it is well-known (cf.\cite{ADN}) that $\nabla u\in L^p_{\loc}(U)$ for $1<p<\infty$; in fact, this is even true when the coefficients are in VMO (cf.\cite{DF}). If the coefficients are Dini-continuous in $U$, then $u$ is known to be continuously differentiable (cf.\cite{HW},\cite{T}).
In the present paper, we find conditions on the coefficients $a_{ij}$, milder than Dini-continuity, under which $u$ must be Lipschitz continuous, or even differentiable, at a given point.

Let us fix an interior point of $U$, which for convenience we shall assume is the origin, $x=0$. Using a change of independent variables, we may assume that $a_{ij}(0)=\de_{ij}$. Suppose that
\begin{equation}
\sup_{|x|=r}|a_{ij}(x)-\de_{ij}|\leq \om(r) \quad\hbox{as $r\to 0$, }
\label{eq:a_ij-delta_ij}
\end{equation}
where $\om(r)$ is a continuous, nondecreasing function for $0\leq r< 1$ satisfying $\om(0)=0$.
We shall not require the Dini condition on $\om$, i.e.\   $r^{-1}\om(r)\in L^1(0,1)$; instead we assume that $\om$ satisfies the {\it square-Dini condition}:
\begin{equation}
\int_0^1 \om^2(r)\,\frac{dr }{r}< \infty.
\label{eq:Sq-Dini}
\end{equation}
However, examples show that additional conditions are required to ensure that a solution is even Lipschitz continuous.

Our additional conditions for regularity are derived from a dynamical system that we shall now describe. Let 
\begin{equation}
 R(r):=\meanint_{S^{n-1}}\left(A(r\theta)-nA(r\theta)\theta\otimes\theta\right)\,ds_\theta,
\label{def:R}
\end{equation}
where the slashed integral denotes mean value, $A=(a_{ij})$, $r=|x|$, $\theta=x/|x|\in S^{n-1}$, $A\theta\otimes\theta$ is the outer product of the vectors $A\theta$ and $\theta$, and $ds$ denotes standard surface measure on $S^{n-1}$. Note that  $\pmb{|}R(r)\pmb{|}\leq c\,\om(r)$, where  we use $\pmb{|}\cdot\pmb{|}$ to denote the matrix norm. Also note that $R$ need not be symmetric. Let us consider the dynamical system
\begin{equation}
 \frac{d\phi}{dt} + R(e^{-t})\,\phi=0 \quad\hbox{for}\ T<t<\infty,
 \label{eq:DynSystm}
 \end{equation}
where $t=-\log r$ and $T$ is sufficiently large. We shall find that the regularity of weak solutions of (\ref{eq:Lu=0}) is determined by the asymptotic behavior as $t\to\infty$ of solutions of (\ref{eq:DynSystm}). We say that  (\ref{eq:DynSystm}) is {\it uniformly stable} as $t\to \infty$ if for every $\e>0$ there exists a $\de=\de(\e)>0$ such that
any solution $\phi$ of  (\ref{eq:DynSystm}) satisfying $|\phi(t_1)|< \de$ for some $t_1>0$ satisfies $|\phi(t)|<\e$ for all $t\geq t_1$ (cf.\ \cite{C}). In addition, we are interested in the condition that {\it every solution of  (\ref{eq:DynSystm}) is  asymptotically constant}, i.e. $\phi(t)\to\phi_\infty$ as $t\to\infty$. These two stability conditions are independent of each other (cf.\ Section 5).  On the other hand, it is easy to see that  $r^{-1}R(r)\in L^1(0,\e)$ implies that (\ref{eq:DynSystm}) is uniformly stable and every solution is asymptotically constant  as $t\to \infty$; in particular, if $\om$ satisfies the Dini condition, then these conditions are met.

 We are now in a position to state the main result of this paper; since we are only concerned with regularity at $x=0$, the coefficients are not required to be continuous elsewhere.
 \begin{theorem} \label{th:main}
Suppose that $a_{ij}$ satisfy (\ref{eq:a_ij-delta_ij}) where $\om$ satisfies (\ref{eq:Sq-Dini}) and that (\ref{eq:DynSystm}) is uniformly stable. Then every weak solution $u\in H_{\loc}^{1,2}(U)$ of (\ref{eq:Lu=0}) is Lipschitz continuous at $x=0$ and satisfies
 \begin{equation}
 |u(x)-u(0)|\leq \frac{c\,|x|}{r}\left(\meanint_{|y|<r}|u(y)|^2\,dy\right)^{1/2}
 \quad\hbox{for}\ |x|<r/2,
 \label{est:LipBnd}
 \end{equation}
where $r$ is sufficiently small. In addition, if every solution of (\ref{eq:DynSystm}) is asymptotically constant, then $u$ is differentiable at $x=0$ and we have
 \begin{equation} \label{eq:del_j_u(0)}
 \del_j u(0)=\lim_{r\to 0}\,\frac{n}{r}\,\meanint_{S^{n-1}} u(r\theta)\,\theta_j\,ds_\theta.
 \end{equation}
 \end{theorem}
 \noindent
 In this theorem and throughout this paper, $c$ denotes a constant whose value may change with the instance but does not depend upon the solution $u$ or the parameter $r$.
 
 \begin{remark} If the $a_{ij}$ are radial functions, then $R(r)\equiv 0$ and we only require (\ref{eq:a_ij-delta_ij}) and (\ref{eq:Sq-Dini}) to conclude that weak solutions are differentiable at $x=0$. Moreover, if $a_{ij}(x)=a^0_{ij}(|x|)+a^1_{ij}(x)$, then the $R$ in (\ref{eq:DynSystm}) is completely determined by $a^1_{ij}$; for example, if the $a^1_{ij}$ are Dini continuous then weak solutions are differentiable  even though $a_{ij}$ need only be square-Dini continuous.
\label{rk:1} \end{remark}
 
We also investigate specific analytic conditions on the coefficients $a_{ij}$ that imply the desired asymptotic properties of (\ref{eq:DynSystm}). 
  Let  us introduce the symmetric matrix ${\mathcal S}=-\frac{1}{2}(R+R^t)$, i.e.\ 
  \begin{equation}
{\mathcal S}(r):=\meanint_{S^{n-1}}\left(\frac{n}{2}\left[A(r\theta)\theta\otimes\theta+\theta\otimes A(r\theta)\theta\right]-A(r\theta)\right)\,ds_\theta,
\label{def:mathcal-R}
\end{equation}
and
\begin{equation}
\mu({\mathcal S})=\hbox{largest eigenvalue of} \ {\mathcal S}.
\label{def:mu}
\end{equation}
In Section 4 we use the theory of dynamical systems to show that if 
 there exist positive constants $\e$ and $K$ so that
 \begin{equation}
 \int_{r_1}^{r_2} \mu\left({\mathcal S}(\rho)\right) \,\frac{d\rho}{\rho}<K \quad\hbox{for all}\ \e> r_2 > r_1 >0,
 \label{eq:Lipschitz-condition}
 \end{equation}
then (\ref{eq:DynSystm}) is uniformly stable. As a consequence, Theorem 1 implies the following:
 
 \begin{corollary}     Suppose that  (\ref{eq:a_ij-delta_ij}),  (\ref{eq:Sq-Dini}), and 
 (\ref{eq:Lipschitz-condition}) are satisfied.
Then every weak solution $u$ of (\ref{eq:Lu=0}) is Lipschitz continuous at $x=0$.
\label{cor:1}
 \end{corollary}

\noindent
What about conditions for differentiability at $x\!=\!0$\,? As already observed, $r^{-1}R(r)\in L^1(0,\e)$ is sufficient, but is there a weaker condition? 
Let us suppose that for $r\in(0,\e)$ the improper integral
  \begin{subequations}\label{Differentiability:joint}
  \begin{equation}\label{Differentiability:1}
  \int_0^r R(\rho)\,\frac{d\rho}{\rho} \quad\hbox{converges (perhaps not absolutely)}.
  \end{equation}
Examples show (see Section 5) that this condition is not sufficient to ensure that (\ref{eq:DynSystm}) is uniformly stable; we shall require an additional condition such as
  \begin{equation}\label{Differentiability:2}
 \frac{R(r)}{r}\int_0^r  R(\rho)\,\frac{d\rho}{\rho} \in L^1(0,\e),
\end{equation}
\end{subequations}
which is also weaker than assuming $R(r)\,r^{-1}\in L^1(0,\e)$. In Section 4 we show that (\ref{Differentiability:1}) and (\ref{Differentiability:2}) together imply not only that (\ref{eq:DynSystm}) is uniformly stable but asymptotically constant.
Consequently, Theorem 1 yields the following:

 \begin{corollary}    Suppose that  (\ref{eq:a_ij-delta_ij}) and  (\ref{eq:Sq-Dini}) are satisfied,
as well as both (\ref{Differentiability:1}) and (\ref{Differentiability:2}).
Then every weak solution $u$ of (\ref{eq:Lu=0}) is differentiable at $x=0$. 
\label{cor:2}
 \end{corollary}
  
  \begin{remark}  Just as (\ref{Differentiability:1}) and (\ref{Differentiability:2})  have replaced the more restrictive $R(r)\,r^{-1}\in L^1(0,\e)$, the assumption (\ref{Differentiability:2}) may be replaced by 
  \begin{subequations}\label{Differentiability':joint}
  \begin{equation}
\int_0^r R(\rho)\left(\int_0^\rho  R(\si)\,\frac{d\si}{\si}\right)\frac{d\rho}{\rho} 
\quad\hbox{converges (perhaps not absolutely)},
\label{Differentiability':1}
\end{equation}
and
  \begin{equation}
\frac{R(r)}{r}\int_0^r R(\rho)\left(\int_0^\rho  R(\si)\,\frac{d\si}{\si}\right)\frac{d\rho}{\rho}\in L^1(0,\e).
\label{Differentiability':2}
\end{equation}
\end{subequations}
This process may be iterated to obtain further refinements.
\label{rk:2} \end{remark}
 
 \begin{remark}  The condition  (\ref{Differentiability:1}) can be expressed as a volume integral (computed in the sense of Cauchy principal value):
  \begin{equation}\label{Differentiability1'}
  \int_{|x|<r} \left(A(x)-n\,\frac{A(x)x}{|x|}\otimes\frac{x}{|x|}\right)\frac{dx}{|x|^n} \quad\hbox{converges for $r\in(0,\e)$}.
  \end{equation}
  This form of the condition is better suited for changes of coordinates, so can be expressed without the simplifying assumption $a_{ij}(0)=\de_{ij}$; however, (\ref{Differentiability:2}) is not so easily handled in this way. In a similar spirit, the following condition\begin{equation}
\int_{|x|<\e}\pmb{\bigr |}A(x)-I\pmb{\bigr |}\frac{dx}{|x|^n}<\infty
\label{|A-I|-integrable}
\end{equation}
is sufficient for Corollary \ref{cor:2} and
easily generalizes to the case $a_{ij}(0)\not=\de_{ij}$; however, it implies $r^{-1}R(r)\in L^1(0,\e)$, so is less general than assuming (\ref{Differentiability:1}) and (\ref{Differentiability:2}).
\label{rk:3} \end{remark}

 \noindent
 Let us consider one more consequence of Theorem 1. In Section 4 we show that
 \begin{equation}
 \int_r^\e \mu({\mathcal S}(\rho))\frac{d\rho}{\rho}\to -\infty \quad\hbox{as}\ r\to 0
 \label{int-Rgoesto-infty}
 \end{equation}
 implies that the null solution of (\ref{eq:DynSystm}) is asymptotically stable. Thus Theorem 1 yields the following:
 
  \begin{corollary}    Suppose that  (\ref{eq:a_ij-delta_ij}),  (\ref{eq:Sq-Dini}), and 
 (\ref{eq:Lipschitz-condition}) are satisfied. Moreover, assume (\ref{int-Rgoesto-infty}).
 Then every weak solution $u$ of (\ref{eq:Lu=0}) is differentiable at $x=0$ and all derivatives are zero: $\del_ju(0)=0$ for $j=1,\dots,n$.
\label{cor:3}
 \end{corollary}
 \noindent
 In Section 5  we discuss an example which illustrates the sharpness of the conditions in Theorem 1 and its corollaries.
  
\smallskip
Now we describe the main ideas of the proof of Theorem 1, which is given in Section 3.
We  write $u$ in the form
\begin{equation}
u(x)=u_0(|x|)+ \vec v(|x|) \cdot\vec x + w(x), 
\label{eq:spectral_decomp}
\end{equation}
where 
\begin{equation}
u_0(r):=\meanint_{S^{n-1}} u(r\theta)\,ds_\theta,
\qquad
v_k(r):=\frac{n}{r}\,\meanint_{S^{n-1}} u(r\theta)\theta_k\,ds_\theta, 
\label{def:uk}
\end{equation}
and $w$ has zero spherical mean and first spherical moments:
\begin{equation}
\meanint_{S^{n-1}}w(r\theta)\,ds_\theta=0=\meanint_{S^{n-1}}w(r\theta)\theta_i\,ds_\theta\quad\hbox{for}\ i=1,\dots,n.
\label{eq:meanint-w}
\end{equation}
We shall find that $\vec v$ satisfies a second-order differential system depending upon $u_0$ and $w$, but it is equivalent to a first-order system that only depends on $w$. Moreover, in this first-order system, 
the behavior of both $\vec v(r)$ and $r\vec v{\,}'(r)$ is controlled  by the asymptotic properties of the solutions to (\ref{eq:DynSystm}). 
To be more specific, we need to assume that $\om(r)$ does not vanish faster than $r$ as $r\to 0$:  
\begin{equation}
 \om(r)\,r^{-1+\kappa} \quad \hbox{is nonincreasing for $r$ near $0$ and some $\kappa>0$.}
\label{om*r(-1+k)}
\end{equation}
Then  the assumption that (\ref{eq:DynSystm}) is uniformly stable ensures not only that $\vec v(r)$ and $r\vec v{\,}'(r)$ are bounded as $r\to 0$, but that  $|u_0(r)-u_0(0)|$ and $|w(x)|$ are both bounded by  $c\,\om(r)\,r$, so we obtain
\begin{equation}
u(x)=u(0)+ \vec v(|x|) \cdot\vec x + O(\om(r)r) \quad\hbox{as}\ r=|x|\to 0, 
\label{eq:Taylor1}
\end{equation}
which confirms that $u$ is Lipschitz at $x=0$. If we also know that all solutions of (\ref{eq:DynSystm}) are asymptotically constant, then $\vec v(r)=\vec v(0)+o(1)$ as $r\to 0$, which shows that $u$ is differentiable at $x=0$.
 
 We observe that the square-Dini condition has been encountered by several other authors in a variety of  contexts. It was used by Stein and Zygmund \cite{SZ} in their investigation of the differentiability of functions, by Fabes, Sroka, and Widman \cite{FSW} in their study of Littlewood-Paley estimates for parabolic equations, and more recently by several authors (cf.\ \cite{FJK}, \cite{CL}, \cite{F}, \cite{KP}) investigating the absolute continuity of elliptic measure and $L^2$ boundary conditions for the Dirichlet problem. We also used the square-Dini condition in \cite{MM1} and \cite{MM2} to study  equations in nondivergence and double divergence form.  
 In addition, we should note that the techniques used in this paper are related to, but independent of,  the asymptotic theory developed in ~\cite{KM}.

Finally, we mention that the techniques and results of this paper apply to weak solutions of more general linear equations than (\ref{eq:Lu=0}):  lower-order terms in $u$ (even with mild singularities in the coefficients) as well as a nonhomogenous right-hand side (with certain integrability conditions) can be treated. However, sharp conditions on singular coefficients in lower-order terms requires additional analysis beyond the results of this paper.

\section{Potential Theory Estimates in $\RR^n$}

We will encounter an equation in the following form
\begin{equation}
-\Lap w=g \qquad\hbox{in}\ \RR^n\backslash\{0\},
\label{eq:Laplace-w=g}
\end{equation}
where $g$ is a distribution and we consider the derivatives in  (\ref{eq:Laplace-w=g}) in the distributional sense. We will encounter certain orthogonality conditions with respect to the spherical mean, so let us summarize these in the following:

\begin{lemma} 
(a) Suppose that $f\in L_{\loc}^1(\RR^n\backslash\{0\})$ and $g$ is bounded with compact support. Then
\begin{subequations} \label{L1:joint}
 \begin{equation}
\int_{\RR^n}\overline{f}(|x|)g(x)\,dx=\int_{\RR^n}{f}(x)\overline g(|x|)\,dx.
\label{L1:a}
\end{equation}
(b) Suppose that $\nabla f\in L_{\loc}^1(\RR^n\backslash\{0\})$. For all $r>0$
 \begin{equation}
\overline{f}(r)=\meanint_{S^{n-1}}  f(r\theta)\,ds_\theta=0 \ \Rightarrow\ 
\meanint_{S^{n-1}}  \, \theta_i\,\del_if(r\theta)\,ds_\theta=0.
\label{L1:b}
\end{equation}
(c) Suppose that $\nabla f\in L_{\loc}^1(\RR^n\backslash\{0\})$.  For all $r>0$ and any $i=1,\dots,n$
 \begin{equation}
\meanint_{S^{n-1}}  \,\theta_i f(r\theta)\,ds_\theta=0 \ \Rightarrow\ 
\meanint_{S^{n-1}}  \,\del_i f(r\theta)\,ds_\theta=0=\meanint_{S^{n-1}}  \,\theta_i\theta_j\del_j f(r\theta)\,ds_\theta.
\label{L1:c}
\end{equation}
\end{subequations}
\label{le:lemma}
\end{lemma}

\noindent{\bf Proof.}
The proof of (a) is trivial. To prove (b), we consider $\phi\in C_0^\infty(0,\infty)$ and compute
\[\begin{aligned}
\biggl\langle \meanint \theta_i\del_i f ds_\theta,\phi\biggr\rangle
&=\int_0^\infty\meanint\,\theta_i\,\del_if(r\theta)ds_\theta\,\phi(r)\,dr
=\frac{1}{|S^{n-1}|}\int_{\RR^n} x_i\del_i f(x)\phi(|x|)|x|^{-n}\,dx\\
&=-\frac{1}{|S^{n-1}|}\int_{\RR^n} nf(x)\phi(|x|)|x|^{-n}\,dx
-\frac{1}{|S^{n-1}|}\int_{\RR^n}x_if(x)[\phi(r)r^{-n}]'|_{r=|x|}\theta_i\,dx\\
&=-\int_0^\infty\left(\meanint_{S^{n-1}}f(r\theta)ds_\theta\right)\phi'(r)\,dr=0
\end{aligned}\]
where $'$ denoted $d/dr$. To prove (c), we again consider $\phi\in C_0^\infty(0,\infty)$ and compute
\[\begin{aligned}
\biggl\langle \meanint \,\del_i f ds,\phi\biggr\rangle
&=\int_0^\infty \meanint_{S^{n-1}}\del_i f(r\theta)\,ds\,\phi(r)\,dr
=\frac{1}{|S^{n-1}|} \int_{\RR^n}\del_if(x)|x|^{1-n}\phi(|x|)\,dx\\
&=-\frac{1}{|S^{n-1}|} \int_{\RR^n} f(x)[r^{1-n}\phi(r)]'|_{r=|x|}\theta_i\,dx\\
&=-\int_0^\infty\int_{S^{n-1}}f(r\theta)\theta_i\,ds\,[r^{1-n}\phi(r)]' r^{n-1}dr=0.
\end{aligned}\]
The proof of the remaining identity in (\ref{L1:b}) is similar.
\hfill$\Box$

\smallskip
Note that (\ref{L1:a}) enables us to define the spherical mean of a distribution. In fact, for $f\in L_{\loc}^1(\RR^n\backslash\{0\})$, let us define 
\begin{equation}
{f(r\theta)}^\perp=f(r\theta)-Pf(r\theta),
\label{def:perp}
\end{equation}
where $Pf$ is defined by
\[
Pf(r\theta)=\meanint_{S^{n-1}} f(r\phi)\,ds_\phi+n\theta_k\meanint_{S^{n-1}}\phi_kf(r\phi)\,ds_\phi.
\]
Using
\[
\meanint_{S^{n-1}} \ \theta_k\theta_\ell\,ds_\theta=\frac{1}{n}\,\de_{k\ell} \quad\hbox{for}\ k, \ell=1,\dots,n,
\]
it is clear that  (for each $r>0$) $P$ is the projection of $f$ onto the functions on $S^{n-1}$ spanned by $1, \theta_1,\dots,\theta_n$. For the same reason as (\ref{L1:a}), we have
\begin{equation}
\int_{\RR^n}P(f)\,g\,dx=\int_{\RR^n}{f}\, P(g)\,dx.
\label{eq:int-P}
\end{equation}
for $g$ bounded with compact support.
This also allows us to define $P$ on distributions. In fact, one particular instance of 
(\ref{eq:Laplace-w=g}) that we are interested in is
\begin{equation}
-\Lap w=[\del_i f_i]^\perp\qquad\hbox{in}\ \RR^n\backslash\{0\},
\label{eq:Laplace-w=df-P(df)}
\end{equation}
where $f_i\in L^1_{\loc}(\RR^n\backslash\{0\})$. We will solve (\ref{eq:Laplace-w=df-P(df)}) by convolution with $\Gamma(|x|)$, the fundamental solution for $-\Lap$, but
we are interested in controlling the growth of the solution near $x=0$ using mean values over annuli. We consider $L^p$-means for any $p\in (1,\infty)$:
 \begin{equation}
M_{p}(w,r)=\left(\meanint_{A_r}|w(x)|^p\,dx\right)^{1/p},
\label{def:Mp}
\end{equation}
where $A_r=\{x:r<|x|<2r\}$ and $w$ may be scalar or vector-valued. To control the growth of the first derivatives of functions we introduce
 \[
M_{1,p}(w,r)=rM_p(\nabla w,r)+M_p(w,r).
\]

\begin{proposition}
Suppose $n\geq 2$, $p\in (1,\infty)$, and $f=(f_i)$ with $f_i\in L_{\loc}^p(\RR^n\backslash\{0\})$ satisfies
\[
\int_{|x|< 1} |x|\,|f(x)|\,dx<\infty, 
\quad\hbox{and}\quad
\int_{|x|> 1} |f(x)|\,|x|^{-n-1}\,dx<\infty.
\]
Then convolution by $\Gamma$ defines a solution $w\in H^{1,p}_{\loc}(\RR^n\backslash\{0\})$  of (\ref{eq:Laplace-w=df-P(df)}) that satisfies $Pw=0$ and
\begin{equation}
M_{1,p}(w,r)\leq c
\left(r^{-n}\int_{0}^r\,M_p(f,\rho)\,\rho^{n}\,d\rho+
r^2\int_{r}^\infty M_p(f,\rho)\,\rho^{-2}\,d\rho
\right).
\label{eq:M1est-w}
\end{equation}
\label{pr:Prop1}
\end{proposition}
 \noindent In (\ref{eq:M1est-w}) and throughout this paper, $c$ denotes a constant; in other instances, the value of $c$ may change line by line without change in notation.

\smallskip
\noindent{\bf Proof of Proposition 1.} We may assume that $f_i\in C_0^1(\RR^n\backslash\{0\})$ since the general case may be handled by an approximation argument. 
Using (\ref{eq:int-P}), let us write 
the solution of  (\ref{eq:Laplace-w=df-P(df)}) as
\[
\begin{aligned}
w(x)=&
\int_{\RR^n}\Gamma(|x-y|)\left(\del_if_i(y)-P( \del_if_i)(y)\right)\,dy\\
&=\int_{\RR^n}\left[\Gamma(|x-y|)-P(\Gamma_x)(|y|,\hat y)\right]\del_if_i(y)\,dy,
\quad\hat y=y/|y|,
\end{aligned}
\]
where $\Gamma_x(y)=\Gamma(|x-y|)$; clearly $Pw=0$. To calculate $P(\Gamma_x)$, we use an expansion of $\Gamma$ in spherical harmonics. Let ${\mathcal H}_k$ denote the spherical harmonics of degree $k$ and let $N(k)=\hbox{dim}{\mathcal H}_k$. For each $k$, choose an orthonormal basis $\{\phi_{k,m}:m=1,\dots,N(k)\}$ for ${\mathcal H}_k$; for $k=1$, note that $\phi_{1,m}(\theta)=\sqrt{n}\,\theta_m$. For notational convenience, let us assume $n\geq 3$; the case $n=2$ is analogous.
For $|x|<|y|$ we can write $\Gamma(|x-y|)$ as a convergent series
\begin{equation}
\Gamma(|x-y|)=\sum_{k=0}^\infty \,\frac{|x|^k}{|y|^{n-2+k}}\sum_{m=1}^{N(k)}\phi_{k,m}\left(\hat{x}\right)\,\phi_{k,m}\left(\hat{y}\right).
\label{eq:Gamma=expansion}
\end{equation}
 Since $\int\phi_{k,m}(\theta)\,ds_\theta=0$ for $k>0$, the spherical mean of $\Gamma$ is given by
\[
\meanint_{S^{n-1}}\Gamma(|x-y|)\,ds_{\hat{y}}=c_n\,|y|^{2-n}=\Gamma(|y|).
\]
We can also use (\ref{eq:Gamma=expansion}) and the orthogonality of the $\phi_{k,\ell}$ to compute
\[
\meanint_{S^{n-1}}\theta_\ell\,\Gamma_x(|y|\theta)\,ds_\theta=\frac{|x|}{|y|}\Gamma(|y|)\hat x_\ell.
\]
Consequently,
\[
P(\Gamma_x)(|y|,\hat y)=\left(1+n\frac{|x|}{|y|}\hat x\cdot\hat y\right)\Gamma(|y|).
\]
By symmetry, it is clear how to modify these projections for $|x|>|y|$, so we obtain the following:
\[
\begin{aligned}
w(x)=&\int_{|y|>|x|} \left[\Gamma(|x-y|)-\Gamma(|y|)-n\,\frac{|x|}{|y|}\,\Gamma(|y|)\,\hat x\cdot\hat y\right]\,\del_i f_i(y)\,dy\\
&+\int_{|y|<|x|} \left[\Gamma(|x-y|)-\Gamma(|x|)-n\,\frac{|y|}{|x|}\,\Gamma(|x|)\,\hat x\cdot\hat y\right]\,\del_i f_i(y)\,dy\\
=&-\int_{|y|>|x|} \frac{\del}{\del y_i}\left[\Gamma(|x-y|)-\Gamma(|y|)-n\,\frac{|x|}{|y|}\,\Gamma(|y|)\,\hat x\cdot\hat y\right]\, f_i(y)\,dy\\
&-\int_{|y|<|x|}  \frac{\del}{\del y_i}\left[\Gamma(|x-y|)-\Gamma(|x|)-n\,\frac{|y|}{|x|}\,\Gamma(|x|)\,\hat x\cdot\hat y\right]\, f_i(y)\,dy,
\end{aligned}
\]
where we have used the divergence theorem (and the fact that $f_i$ is supported in $\RR^n\backslash\{0\}$).

If we assume $r<|x|<2r$ and introduce the annulus $\tilde A_r=\{x:r/2<|x|<4r\}$, then we can split up the integrals as follows:
\[
\begin{aligned}
w(x)=&-\int_{\tilde A_r} \frac{\del}{\del x_i}\Gamma(|x-y|)f_i(y)\,dy+
\int_{r/2<|y|<|x|}\frac{\del}{\del y_i}\left(\Gamma(|x|)+n\frac{|y|}{|x|}\,\Gamma(|x|)\,\hat x\cdot\hat y\right)f_i(y)\,dy\\
& +\int_{|x|<|y|<4r}\frac{\del}{\del y_i}\left(\Gamma(|y|)+n\frac{|x|}{|y|}\,\Gamma(|y|)\,\hat x\cdot\hat y\right)f_i(y)\,dy\\
&-\int_{|y|<r/2} \frac{\del}{\del y_i}\left[\Gamma(|x-y|)-\Gamma(|x|)-n\,\frac{|y|}{|x|}\,\Gamma(|x|)\,\hat x\cdot\hat y\right]\, f_i(y)\,dy\\
&-\int_{|y|>4r}  \frac{\del}{\del y_i}\left[\Gamma(|x-y|)-\Gamma(|y|)-n\,\frac{|x|}{|y|}\,\Gamma(|y|)\,\hat x\cdot\hat y\right]\, f_i(y)\,dy.
\end{aligned}
\]
Using (\ref{eq:Gamma=expansion}) we can estimate the last two integrals:
\[
\begin{aligned}
\left|\int_{|y|<r/2} \frac{\del}{\del y_i}\left[\Gamma(|x-y|)-\Gamma(|x|)-n\,\frac{|y|}{|x|}\,\Gamma(|x|)\,\hat x\cdot\hat y\right]\, f_i(y)\,dy\right|
&\leq c\int_{|y|<r/2}|x|^{-n}|yf(y)|\,dy\\
\leq \ & c\,r^{-n}\int_{|y|<|x|}|y|\,|f(y)|\,dy\\
\left|\int_{|y|>4r}  \frac{\del}{\del y_i}\left[\Gamma(|x-y|)-\Gamma(|y|)-n\,\frac{|x|}{|y|}\,\Gamma(|y|)\,\hat x\cdot\hat y\right]\, f_i(y)\,dy\right|&\leq c\int_{|y|>4r} |x|^2|y|^{-n-1}|f(y)|\,dy\\
\leq \ &c\,r^2\int_{|y|>|x|} |y|^{-n-1}|f(y)|\,dy.
\end{aligned}
\]
We can also easily estimate
\[
\begin{aligned}
\left|\int_{r/2<|y|<|x|}\frac{\del}{\del y_i}\left(\Gamma(|x|)+n\frac{|y|}{|x|}\,\Gamma(|x|)\,\hat x\cdot\hat y\right)f_i(y)\,dy\right|
&\leq c\int_{r/2<|y|<|x|}|x|^{1-n}|f(y)|\,dy\\
\leq c\,r^{-n}\int_{r/2<|y|<|x|}|y|\,|f(y)|\,dy&\leq  c\,r^{-n}\int_{|y|<|x|}|y|\,|f(y)|\,dy
\end{aligned}
\]
and
\[
\begin{aligned}
\left|\int_{|x|<|y|<4r}\frac{\del}{\del y_i}\left(\Gamma(|y|)+n\frac{|x|}{|y|}\,\Gamma(|y|)\,\hat x\cdot\hat y\right)f_i(y)\,dy\right|
\leq \ & c\, r^2 \int_{|x|<|y|<4r} |y|^{-n-1}|f(y)|\,dy\\
\leq \ & c\, r^2 \int_{|x|<|y|} |y|^{-n-1}|f(y)|\,dy.
\end{aligned}
\]
We conclude that
\[
\left|w(x)-\int_{\tilde A_r} \frac{\del}{\del x_i}\Gamma(|x-y|)f_i(y)\,dy\right| 
\leq
c\left(r^{-n}\int_{|y|<|x|}|y|\,|f(y)|\,dy+r^2\int_{|y|>|x|}|y|^{-n-1}|f(y)|\,dy
\right).
\]
Similarly, we can show that
\[
\left|\frac{\del w}{\del x_j}(x)-\int_{\tilde A_r}\frac{\del^2}{\del x_j\del x_i}\Gamma(|x-y|)f_i(y)\,dy\right|
 \]
 is bounded by
\[
c\left( r^{-n-1}\int_{|y|<|x|}|y|\,|f(y)|\,dy+r\int_{|y|>|x|}|y|^{-n-1}|f(y)|\,dy\right).
\]
Now we use Stein's inequality \cite{St} to conclude
\[
\left\|\int_{\tilde A_r}\frac{\del}{\del x_i}\Gamma(|x-y|)f_i(y)\,dy\right\|_{L^p(A_r)}
\leq c\,r\,\|f\|_{L^p(\tilde A_r)}
\]
and the $L^p$-boundedness of singular integral operators to conclude
\[
\left\|\int_{\tilde A_r}\frac{\del^2}{\del x_j\del x_i}\Gamma(|x-y|)f_i(y)\,dy\right\|_{L^p(A_r)}
\leq c\,\|f\|_{L^p(\tilde A_r)}.
\]
Putting this all together, we obtain
\begin{equation}
M_{1,p}(w,r)\leq c\,\left(r\tilde M_p(f,r)+r^{-n}\int_{|y|<r}|y|\,|f(y)|\,dy+r^2\int_{|y|>r}|f(y)||y|^{-n-1}\,dy
\right)
\label{eq:pre-M1est-mu}
\end{equation}
where $\tilde M_p$ denotes the mean value over $\tilde A_r$ instead of $A_r$.

The integrals in (\ref{eq:pre-M1est-mu}) can be estimated in terms of $M_p$ and combined with the $\tilde M_p$ term. In fact it is elementary (cf.\ \cite{MM1}) to establish
\[
\int_{|y|<r}|y|\,|f(y)|\,dy\leq c\int_0^r M_p(f,\rho)\rho^{n}\,d\rho
\]
\[
\int_{|y|>r}|f(y)||y|^{-n-1}\,dy\leq c\int_r^\infty M_p(f,\rho)\rho^{-2}\,d\rho.
\]
In addition, it is easy to see that
\[
\tilde M_p(f,r)\leq 
\,c\int_{r/4}^{4r} M_p(f,\rho)\rho^{-1}d\rho
\leq c\left[r^{-n-1}\int_{r/4}^r M_p(f,\rho)\rho^{n}\,d\rho+r\int_r^{4r} M_p(f,\rho)\rho^{-2}d\rho\right].
\]
Using these inequalities, it is clear that (\ref{eq:pre-M1est-mu}) implies
(\ref{eq:M1est-w}). 
\hfill $\Box$

\medskip
Another instance of 
(\ref{eq:Laplace-w=g}) that we are interested in is
\begin{equation}
-\Lap w=[f]^\perp\qquad\hbox{in}\ \RR^n\backslash\{0\},
\label{eq:Laplace-w=f-P(f)}
\end{equation}
where $f\in L^1_{\loc}(\RR^n\backslash\{0\})$. To control the growth of the second derivatives of functions we use
 \[
M_{2,p}(w,r)=r^2M_p(D^2 w,r)+M_{1,p}(w,r),
\]
where $D^2w$ denotes the Hessian matrix of all second-order derivatives of $w$. The proof of the following is analogous to that of Proposition 1.

\begin{proposition}
Suppose $n\geq 2$, $p\in (1,\infty)$, and $f\in L_{\loc}^p(\RR^n\backslash\{0\})$ satisfies
\[
\int_{|x|< 1} |x|^2\,|f(x)|\,dx<\infty, 
\quad\hbox{and}\quad
\int_{|x|> 1} |f(x)|\,|x|^{-n}\,dx<\infty.
\]
Then convolution by $\Gamma$ defines a solution $w\in H^{2,p}_{\loc}(\RR^n\backslash\{0\})$  of (\ref{eq:Laplace-w=f-P(f)}) that satisfies
\begin{equation}
M_{2,p}(w,r)\leq c
\left(r^{-n}\int_{0}^r\,M_p(f,\rho)\,\rho^{n+1}\,d\rho+
r^2\int_{r}^\infty M_p(f,\rho)\,\rho^{-1}\,d\rho
\right).
\label{eq:M2p-w}
\end{equation}
\label{pr:Prop2}
\end{proposition}

\section{Stability Properties of Dynamical Systems}

Let us now consider the result that we require from perturbation theory for systems of ODEs on $T<t<\infty$. Without loss of generality, we assume $T=0$.
First, let us introduce a positive, nonincreasing continuous function $\e(t)$ satisfying
\begin{equation}
\int_0^\infty \e^2(t)\,dt< \infty.
\label{eq:epsilon-Sq-Int}
\end{equation}
Now consider the $2n\times 2n$ system on $(0,\infty)$
\begin{subequations} \label{DS:joint}
\begin{equation} \label{DS:system}
\frac{d}{dt}
\begin{pmatrix}
\phi \\ \psi
\end{pmatrix} + \begin{pmatrix}
0 & 0 \\ 0 & -nI
\end{pmatrix}
\begin{pmatrix}
\phi \\ \psi
\end{pmatrix} + {\mathcal R}(t)
\begin{pmatrix}
\phi \\ \psi
\end{pmatrix} =  g(t),
\end{equation}
where i) ${\mathcal R}$ is a $2n\times 2n$ matrix of the form
\begin{equation}\label{DS:R}
{\mathcal R}(t)=
\begin{pmatrix}
R_{1}(t) & R_{2}(t) \\ R_{3}(t) & R_{4}(t)
\end{pmatrix}
\quad
\hbox{with}\ \pmb{\bigr |}R_j(t)\pmb{\bigr |}\leq \e(t)\ \hbox{on}\ 0<t<\infty,
\end{equation}
and ii)  $g=(g_1,g_2)$ with $g_1 \in L^1(0,\infty)$ and there exists $\de>0$ so that for any choice of $\alpha\in [n-\de,n)$ there is a constant $c_\alpha$ so that
\begin{equation}
e^{\alpha t}\int_t^\infty|g_2(s)|e^{-\alpha s}\,ds \leq c_\alpha\,\e(t) \quad\hbox{for}\ 0<t<\infty.
\label{DS:g2}
\end{equation}
\end{subequations}
(With regard to convergence at infinity, (\ref{DS:g2}) is weaker than assuming $g_2\in L^1(0,\infty)$.)
In addition, we assume asymptotic conditions on the solutions of
\begin{equation}
\frac{d\phi}{dt}
 +R_{1}\phi=0 \quad\hbox{for} \ t>0,
\label{eq:R1-unifstable}
\end{equation}
and that $\psi$ satisfies the ``finite-energy condition''
  \begin{equation}
\int_0^\infty\left(|\psi|^2+| {\psi}_t|^2\right)e^{-nt}\,dt<\infty.
 \label{eq:energycondition-psi}
 \end{equation}

\begin{proposition}
Suppose that  (\ref{eq:R1-unifstable}) is uniformly stable. Then all solutions 
$(\phi,\psi)$ of (\ref{DS:joint}) that satisfy (\ref{eq:energycondition-psi})
 remain bounded as $t\to \infty$,  and $\psi(t)\to 0$. In fact, for $\alpha=n-\de$ with $\de>0$ sufficiently small, we have the estimates
 \begin{subequations}\label{sup-bound:joint}
 \begin{equation}
 \sup_{0<t<\infty}|\phi(t)|\leq c\,(c_\alpha+|\phi(0)|+\|g_1\|_1),
 \label{sup-bound:a}
 \end{equation}
 and
 \begin{equation}
| \psi(t)|\leq c\,\e(t)\,(c_\alpha+\sup_{t<\tau<\infty}|\phi(\tau)|\,).
 \label{sup-bound:b}
 \end{equation}
\end{subequations}
 In addition, if all solutions of (\ref{eq:R1-unifstable}) are asymptotically constant as $t\to \infty$, then  the solution $(\phi,\psi)$ of (\ref{DS:joint}) also has a limit:
 \begin{equation}
(\phi(t),\psi(t))\to (\phi_\infty,0) \quad\hbox{as}\ t\to\infty.
 \label{eq:phi-asym}
 \end{equation}
\label{pr:Prop3}
\end{proposition}

\noindent{\bf Proof.}
Let us simplify notation by  denoting $d/dt$ by dot: $d\phi/dt=\dot\phi$. Therefore, we want to study solutions of
\begin{equation}
\begin{aligned}
\dot{\phi}+R_1\phi+R_2\psi&=g\cr
\dot{\psi}-n\psi+R_3\phi+R_4\psi&=h,
\end{aligned}
\label{eq:reduced_system}
 \end{equation}
 when $g\in L^1(0,\infty)$ and $h$ satisfies the condition on $g_2$ in (\ref{DS:g2}) for a certain value of $\alpha$ that will be specified below.
 Let $\Phi$ denote the fundamental matrix for $\dot\phi+R_1\phi=0$ on $t>0$, i.e.\ 
 \[
 \dot{\Phi}+R_1\Phi=0, \quad \Phi(0)=I.
 \]
 The assumption that (\ref{eq:R1-unifstable}) is uniformly stable is equivalent (cf.\ \cite{C}) to 
   \begin{equation}
  \pmb{\bigr |}\Phi(t)\Phi^{-1}(s)\pmb{\bigr |}\leq K\ \hbox{for}\ t>s>0,
   \label{eq:Phi-unifstab}
   \end{equation}
where $K$ is a constant. Next let $\Psi$ denote the fundamental matrix for $\dot\psi+R_4\psi=0$ on $t>0$, i.e.\ 
 \[
 \dot{\Psi}+R_4\Psi=0, \quad \Psi(0)=I.
 \]
Since $\e(t)\to 0$ as $t\to\infty$, for fixed $0<\de<1$  we can find $t_1$ so that 
\begin{equation}
\e(t)<\de\ \hbox{for}\  t\geq t_1, \ \hbox{and}\  \int_{t_1}^\infty \e^2(t)\,dt<\de.
\label{eq:epsilon<delta}
\end{equation}
Without loss of generality, we can assume $t_1= 0$. Using Gronwall's inequality, we can show     \begin{equation}
  \pmb{\bigr |}\Psi(t)\Psi^{-1}(s)\pmb{\bigr |}\leq e^{\de|t-s|} \ \hbox{for}\ t,s>0.
     \label{eq:Psi-stable}
   \end{equation}
However, we also need a lower bound on $\Psi(t)$. To derive this, let $\psi(t)=\Psi(t)\psi_0$ and $p(t)=|\psi(t)|^2$. Then\begin{align} \notag
\dot p = &2\psi_1\dot\psi_1+\dots+2\psi_n\dot\psi_n\\ \notag
= & -2r_{11}\psi_1^2-2r_{12}\psi_1\psi_2-\dots -2r_{1n}\psi_1\psi_n\\ \notag
& -2r_{21}\psi_2\psi_1-2r_{22}\psi_2^2-\dots -2r_{2n}\psi_2\psi_n\\ \notag
&\vdots\\ \notag
& -2r_{n1}\psi_n\psi_1-\dots -2r_{nn}\psi_n^2\\ \notag
\geq & -2n\de |\psi|^2 = -2n\de p.
\end{align}
   Integration yields $p(t)\geq p_0\, e^{-2n\de t}$; in other words
\begin{equation}
\pmb{|}\Psi(t)\pmb{|}\geq e^{-n\de t}.
\label{eq:Psi-lowerbound}
\end{equation}
   
Now we can use the ``variation of parameters'' formula to conclude from the first equation in (\ref{eq:reduced_system}) that
  \begin{equation}
\phi(t)=\Phi(t)\left[\phi(0)+\int_0^t \Phi^{-1}(\tau)[g(\tau)-R_2(\tau)\psi(\tau)]\,d\tau\right],
\label{eq:phi=}
\end{equation}
and from the second equation in (\ref{eq:reduced_system}) that
\begin{equation}
\psi(t)=e^{nt}\Psi(t)\left[\psi(0)+\int_0^t\Psi^{-1}(\tau)[h(\tau)-R_3(\tau)\phi(\tau)]e^{-n\tau}\,d\tau\right].
\label{eq:psi=}
\end{equation}
In order to have (\ref{eq:energycondition-psi}), we see from (\ref{eq:Psi-lowerbound}) that we must have
\[
\psi(0)=-\int_0^\infty \Psi^{-1}(\tau)[h(\tau)-R_3(\tau)\phi(\tau)]e^{-n\tau}\,d\tau,
\]
and consequently (\ref{eq:psi=}) can be rewritten as
\begin{equation}
\psi(t)=e^{nt}\Psi(t)\int_t^\infty \Psi^{-1}(\tau)[R_3(\tau)\phi(\tau)-h(\tau)]e^{-n\tau}\,d\tau.
\label{eq:psi=2}
\end{equation}
If we plug (\ref{eq:psi=2}) into (\ref{eq:phi=}), we obtain 
\begin{subequations}\label{phi-equation:joint}
\begin{equation}
\phi(t)+S\phi(t)=\xi(t)=\xi_0(t)
+\xi_1(t)+\xi_2(t),
\label{phi-equation:1}
\end{equation}
where
\begin{equation}
\begin{aligned}
S\phi(t)&=-\Phi(t)\int_0^t\Phi^{-1}(\tau)R_2(\tau)e^{n\tau}\Psi(\tau)\int_\tau^\infty\Psi^{-1}(\si)R_3(\si)\phi(\si)e^{-n\si}\,d\si\, d\tau,\\
\xi_0(t)&= \Phi(t)\phi(0), \qquad \xi_1(t)=\Phi(t)\int_0^t\Phi^{-1}(\tau)g(\tau)\,d\tau,\\
\xi_2(t)&=\Phi(t)\int_0^t\Phi^{-1}(\tau)R_2(\tau)e^{n\tau}\Psi(\tau)\int_\tau^\infty\Psi^{-1}(\si)h(\si)e^{-n\si}\,d\si\, d\tau.
\end{aligned}
\label{phi-equation:2}
\end{equation}
\end{subequations}
We want to use (\ref{phi-equation:joint}) to conclude that $\phi$ is bounded.

Let $X=C[0,\infty)$ with $\|\phi\|_X:=\sup_{0<t<\infty}|\phi(t)|<\infty$. Notice that (\ref{eq:Phi-unifstab}) and $g\in L^1(0,\infty)$ imply that $\xi_0,\xi_1\in X$. To show $\xi_2\in X$, let us use (\ref{eq:Phi-unifstab}), (\ref{eq:Psi-stable}), and (\ref{DS:g2}):
\[
\begin{aligned}
|\xi_2(t)|&\leq K\int_0^t\e(\tau)e^{(n-\de)\tau}\int_\tau^\infty |h(\si)|e^{(\de-n)\si}\,d\si \,d\tau\\
&\leq K\,\int_0^t\e(\tau)\int_{\tau}^\infty|h(\si)|\,d\si\,d\tau
\leq K\, c_\alpha\,\int_0^t\e^2(\tau) \, d\tau\leq K\,c_\alpha\,\de,
\end{aligned}
\]
although we do not care if this is small. Now let us show that $S:X\to X$ with $\|S\|<1$ if $\de$ is small.
In fact, assume
$\|\phi\|_X\leq 1$. Then
\[
\begin{aligned}
|S\phi(t)|&\leq K\int_0^t \e(\tau)e^{(n-\de)\tau}\int_\tau^\infty e^{(\de-n)\si}\e(\si)|\phi(\si)|\,d\si\,d\tau\\
&\leq K\int_0^t \e(\tau)e^{(n-\de)\tau}\int_\tau^\infty e^{(\de-n)\si}\e(\si)\,d\si\,d\tau \leq \frac{K}{n-\de}\int_0^\infty\e^2(\tau)\,d\tau<\frac{K\de}{n-\de}\, ,
\end{aligned}
\]
so $\|S\|<1$ if $\de$ is sufficiently small. We conclude that $\phi=(1+S)^{-1}\xi\in X$, i.e.\ $\phi$ is bounded.

Now if we apply (\ref{eq:Psi-stable}) to (\ref{eq:psi=2}), we find that
\[
\begin{aligned}
|\psi(t)|&\leq c\, e^{(n-\de)t}\int_t^\infty e^{(\de-n)\tau}\e(\tau)|\phi(\tau)|\,d\tau+e^{(n-\de)t}\int_t^\infty|h(\tau)|e^{(\de-n)\tau}\,d\tau\\
&\leq c\,\e(t) \left(c_\alpha+\sup_{t<\tau<\infty}|\phi(\tau)|\right),
\end{aligned}
\]
which is (\ref{sup-bound:b}). We can then combine this with (\ref{eq:phi=}) and (\ref{eq:epsilon<delta}) to conclude
\[
|\phi(t)|\leq c\,\left(|\phi(0)|+\|g\|_1\right)+c\,\de\left(c_\alpha+\|\phi\|_X\right).
\]
 Taking $\de$ sufficiently small, we can conclude (\ref{sup-bound:a}).

If all solutions of (\ref{eq:R1-unifstable}) are asymptotically constant as $t\to \infty$, then $\Phi(\infty)=\lim_{t\to\infty}\Phi(t)$ exists and 
from (\ref{eq:phi=}) we find
\[
\phi_\infty=\Phi(\infty)\left(\phi(0)+\int_0^\infty\Phi^{-1}(\tau)[g(\tau)-R_2(\tau)\psi(\tau)]d\tau\right).
\]
To show $\phi(t)\to\phi_\infty$ we estimate three terms:
$
\left|\left(\Phi(t)-\Phi(\infty)\right)\phi(0)\right |\leq  \pmb{|}\Phi(t)-\Phi(\infty)\pmb{|} |\phi(0)|,
$
\[
\left|\Phi(t)\int_0^t\Phi^{-1}(\tau)g(\tau)\,d\tau-\Phi(\infty)\int_0^\infty\Phi^{-1}(\tau)g(\tau)\,d\tau\right|\leq \, c\,\left( \pmb{|}\Phi(t)-\Phi(\infty)\pmb{|}\,\|g\|_1 +\int_t^\infty|g(\tau)|\,d\tau\right)
\]
\[
\begin{aligned}
\left|\Phi(t)\int_0^t\Phi^{-1}(\tau)R_2(\tau)\psi(\tau)d\tau-\Phi(\infty)\int_0^\infty\Phi^{-1}(\tau)R_2(\tau)\psi(\tau)d\tau\right|&\\
\leq c\,(1+\|\phi\|_X)&\left[\pmb{|}\Phi(t)-\Phi(\infty)\pmb{|}+\int_t^\infty\e^2(\tau)\,d\tau \right].
\end{aligned}
\]
This confirms (\ref{eq:phi-asym}). 
\hfill $\Box$

\medskip
It will also be useful to know that the uniform stability is not affected by perturbation with terms bounded by $\e^2(t)$.

\begin{proposition}
If $\widetilde R(t)-R(t)\in L^1(0,\infty)$, then (\ref{eq:DynSystm})  is uniformly stable
 if and only if the same is true of
\[
\frac{d\phi}{dt}+\widetilde R\phi=0.
\]
The same equivalency applies if the property that all solutions are asymptotically constant is added to uniform stability.
\label{pr:Prop4}
\end{proposition}

\noindent{\bf Proof.}
Let $S=\widetilde R-R$, $\Phi(t)$ denote the fundamental solution for $\frac{d\phi}{dt}+R\phi=0$ on $t>0$, and $\widetilde\Phi$ denote the same for 
$\frac{d\phi}{dt}+\widetilde R\phi=0$. According to \cite{C}, we know that (\ref{eq:Phi-unifstab}) holds
and we want to prove that a similar bound holds for $\widetilde\Phi$. But we can solve
\[
\frac{d\phi}{dt}+R\phi=-S\phi, \quad \phi(s)=\xi
\]
by variation of parameters to obtain
\[
\phi(t)=\Phi(t)\Phi^{-1}(s)\xi-\Phi(t)\int_s^t\Phi^{-1}(\si)S(\si)\phi(\si)\,d\si,
\]
Applying Gronwall's lemma and (\ref{eq:Phi-unifstab}), we obtain
\[
|\phi(t)|\leq K\,\exp\left(k\int_s^t\pmb{|}S(\si)\pmb{|}\,d\si\right) |\xi| \leq K\,\exp\left(k\int_0^\infty\pmb{|}S(\si)\pmb{|}\,d\si\right)\,|\xi|.
\]
But $\phi(t)=\widetilde\Phi(t)\widetilde\Phi^{-1}(s)\xi$, so this last estimate shows $\frac{d\phi}{dt}+\widetilde R\phi=0$ is uniformly stable. The additional property that solutions are asymptotically constant also follows from the variation of parameters formula. \hfill $\Box$

\section{Proof of Theorem 1}

Recall that $\om(r)$ is a continuous, nondecreasing function for $0\leq r< 1$ satisfying $\om(0)=0$, (\ref{eq:Sq-Dini}),
and (\ref{om*r(-1+k)}). Since our result is local, we focus our attention on a very small ball centered at $x=0$; by rescaling, 
we may assume this is $B_1(0)$, the unit ball centered at $x=0$, and that 
(\ref{eq:a_ij-delta_ij}) holds for $0<r<1$; in fact, given any small $\de>0$, we can similarly assume that
\begin{equation}
\int_0^1 \frac{\om^2(r)dr}{r}<\de \quad\hbox{and}\quad \om(1)\leq \de.
\label{eq:int-omega^2<delta}
\end{equation}
 If the $a_{ij}$ are continuous in $U$, then we know by the results of \cite{ADN} that $\nabla u\in L_{\loc}^p(U)$ for all $p\in (1,\infty)$.
More generally, under the small oscillation assumption (\ref{eq:a_ij-delta_ij}), we can fix any $p\in (1,\infty)$ and conclude that $\nabla u\in L_{\loc}^p(U)$ provided $\de=\de(p)$ in (\ref{eq:int-omega^2<delta}) is sufficiently small; cf.\ Corollary 6.2 in \cite{MMS}. Henceforth, we fix
$p>n$ and assume that $\de(p)$ has been chosen small enough that $\nabla u\in L^p(B_1(0))$.

For our analysis, it is advantageous to extend the problem to all of $\RR^n$, so let us redefine and extend the $a_{ij}$ outside of ${B_1(0)}$ by
\begin{equation}
a_{ij}(x)=\de_{ij} \quad\hbox{ for}\ |x|\geq 1. 
\label{eq:aij_for_r>1}
\end{equation}
If we extend $\om(r)$ to be $\om(1)$ for $r>1$, we see that $\om(r)$ is still nondecreasing and $\om(r)\,r^{-1+\kappa}$ is still nonincreasing; in particular, we easily see that
\begin{equation}
M_p(\om f, r)\leq c\,\om(r)\,M_p(f,r)\quad\hbox{for}\ 0<r<\infty,
\label{eq:Mp(omega_f)}
\end{equation}
where $f\in H_{\loc}^p(\RR^n\backslash\{0\})$ and $c$ is independent of $r$.
Now let us introduce a smooth cutoff function $\chi(r)$ which is $1$ for $0\leq r\leq 1/4$ and $0$ for $r\geq 1/2$. Given a weak solution $u$ of (\ref{eq:Lu=0}), we can write
\[
\del_i(a_{ij}\del_j (\chi u))=\del_i(a_{ij}\chi'\theta_j u)+ \chi'\theta_ia_{ij}\del_j u.
\]
Since we are interested in the behavior near $x=0$ where $\chi(|x|) u(x)$ and $u(x)$ agree, we can assume that $u$ is supported on 
$|x|<1/2$ and satisfies
\begin{subequations}\label{Lu=f:joint}
\begin{equation}\label{Lu=f:a}
\del_i(a_{ij}(x)\del_j u)=\del_i f_i +f_0 \qquad\hbox{in}\  \RR^n,
\end{equation}
where $f_i,\,f_0\in L^p(\RR^n)$ are both supported in $1/4\leq|x|\leq 1/2$, and (using (\ref{def:Lu=0}) with $\eta=\chi$)
\begin{equation}\label{Lu=f:b}
\int_{\RR^n} f_0(x)\,dx=0;
\end{equation}
\end{subequations}
 when convenient, we let $\vec f=(f_1,\dots,f_n)$.
Of course, we now must replace (\ref{def:Lu=0}) with the following:
\begin{equation}
\int_{\RR^n} a_{ij}(x)\,\del_j u\,\del_i\eta\,dx=\int_{\RR^n} f_i\del_i\eta\,dx-\int_{\RR^n} f_0\eta\,dx,
\label{def:Lu=f}
\end{equation}
for all $\eta\in C_0^\infty(\RR^n)$. 

Recall the decomposition $u(x)=u_0(r)+\vec v(r)\cdot\vec x+w(x)$ given in (\ref{eq:spectral_decomp}). Orthogonality properties show that
\begin{equation}
\nabla u\in L^2(B_1(0)) \Rightarrow \int_0^1 \left( (u_0')^2+|\vec v|^2+r^2|\vec v{\,}'|^2\right)r^{n-1}\,dr <\infty\ \hbox{and}\  \nabla w\in L^2(B_1(0)).
\label{eq:u-energy}
\end{equation}
We want to show that $|u_0(r)-u_0(0)|\,r^{-1}\to 0$ and $|w(x)|\,|x|^{-1}\to 0$ as $r=|x|\to 0$, so that the differentiability of $u$ is determined by the behavior of $\vec v$ as $r\to 0$.
The strategy is to show that  $\vec v$ satisfies a system of ODEs  which depends upon $w$ and that $w$ satisfies a PDE 
which depends upon $\vec v$. We obtain the system of ODEs by plugging $u(x)=u_0(r)+\vec v(r)\cdot\vec x+w(x)$ into (\ref{def:Lu=f}) and choosing special $\eta\in C_0^\infty(\RR^n)$. For example, taking $\eta=\eta(r)\in C_0^\infty[0,\infty)$, we obtain 
\begin{subequations}\label{ODE1:joint}
\begin{equation}
\int_0^\infty\left(\alpha\, u_0'+r \vec\beta\cdot\vec v{\,}'+\vec\gamma\cdot\vec v
+p[\nabla w]\right)\eta'\,r^{n-1}\,dr=\int_0^\infty\left(\widetilde f{\,}\eta'-\overline{f_0}\eta\right)\,r^{n-1}\,dr,
\label{ODE1:a}
\end{equation}
where 
\begin{equation}\label{ODE1:b}
\begin{aligned}
\alpha(r)=\meanint_{S^{n-1}} a_{ij}(r\theta)\theta_i\theta_j\,ds_\theta,&
\qquad \beta_k(r)=\meanint_{S^{n-1}} a_{ij}(r\theta)\theta_i\theta_j\theta_k\,ds_\theta,
  \\ 
\gamma_j(r)=\meanint_{S^{n-1}} a_{ij}(r\theta)\theta_i\,ds_\theta,& \qquad
p[\nabla w](r)=\meanint_{S^{n-1}}a_{ij}(r\theta)\,\del_jw(r\theta)\,\theta_i\,ds_\theta,\\
\widetilde f(r)=\meanint_{S^{n-1}} f_i(r\theta)\theta_i\,ds_\theta,& \qquad \overline{f_0}(r)=\meanint_{S^{n-1}}f_0(r\theta)\,ds_\theta. 
\end{aligned}
\end{equation}
Using (\ref{eq:a_ij-delta_ij}) and Lemma 1, we see that these terms satisfy 
\begin{equation}
\begin{aligned}
|\alpha(r)-1|,\ & |\vec\beta(r)|,\ |\vec \gamma(r)|\, \leq \om(r)\ \hbox{for $0<r<1$}, \\
|p[\nabla w](r)|& \leq \om(r)\,\meanint_{S^{n-1}} |\nabla w|\,ds \ \hbox{for $0<r<1$}, \\
\widetilde f(r)=&\overline{f_0}(r)=0  \ \hbox{for $0<r<1/4$}, \\
\end{aligned}
\label{ODE1:c}
\end{equation}
\end{subequations}
whereas for $r>1$ we have $\alpha(r)=1$ and all the other terms vanish.
We can integrate (\ref{ODE1:a}) to find
\begin{subequations}\label{ODE1-int:joint}
\begin{equation}
\alpha(r)u_0'(r)+r \vec\beta(r)\cdot\vec v{\,}'(r)+\vec\gamma(r)\cdot\vec v(r)
+p[\nabla w](r)=\vartheta (r),
\label{ODE1-int:a}
\end{equation}
where
\begin{equation}
\vartheta (r)=\widetilde f(r)+r^{1-n}\int_0^r \overline{f_0}(\rho)\,\rho^{n-1}\,d\rho
\label{ODE1-int:b}
\end{equation}
\end{subequations}
has support in $1/4\leq r\leq 1/2$. (Note that  $\vartheta (r)=0$ for $r<1/4$ follows from (\ref{ODE1:a}) with supp $\eta\subset [0,1/4)$ and $\eta(0)=1$, whereas $\vartheta (r)=0$ for $r>1/2$ is a consequence of (\ref{Lu=f:b}).)

Similarly, we can let $\eta=\eta(r)x_\ell$ in  (\ref{def:Lu=f}) and obtain a 2nd-order linear system of ODEs. We can use (\ref{ODE1-int:a}) to eliminate $u_0$ and then reduce the 2nd-order system for $\vec v$ to a 1st-order system for $(\vec v,\vec v_r)$; of course, these systems also depend on $w$. This 1st-order system is simplified by changing independent variables to $t=-\log r$, so we introduce
\begin{equation}
\e(t)=\om(e^{-t}),
\label{def:epsilon}
\end{equation}
which satisfies (\ref{eq:epsilon-Sq-Int}) by (\ref{eq:Sq-Dini}).  In the Appendix, we show that the 1st-order system for $(\vec v,\vec v_t)$ may be converted to the form  (\ref{DS:joint}):
 \begin{subequations}\label{ourODEsystem:joint}
 \begin{equation}
\frac{d}{dt}
\begin{pmatrix}
\phi \\ \psi
\end{pmatrix}
+
\begin{pmatrix}
0 & 0 \\ 0 & -nI
\end{pmatrix}
\begin{pmatrix}
\phi \\ \psi
\end{pmatrix}
+
{\mathcal R}(t)
\begin{pmatrix}
\phi \\ \psi
\end{pmatrix}
=  
g(t,\nabla w)+h(t),
 \label{ourODEsystem:a}
\end{equation}
where ${\mathcal R}\equiv 0$ for $t<0$, but for $t>0$ it is of the form (\ref{DS:R}) with
\begin{equation}\label{ourODEsystem:b}
R_1(t)\approx \meanint_{S^{n-1}}(A-n\,A\,\theta\otimes\theta)\,ds_\theta
\quad\hbox{as}\ t\to\infty,
\end{equation}
where $\approx$ indicates that the difference is bounded by $c\,\e^2(t)$;
 the term  $g(t,\nabla w)\equiv 0$ for $t<0$ but satisfies
\begin{equation}
|g(t,\nabla w)|\leq c\, \e(t)\,\meanint_{S^{n-1}}  |\nabla w|\,ds\quad\hbox{for}\ t>0,
 \label{ourODEsystem:c}
\end{equation}
and the term $h$ is in $L^p$ with support in $\log 2 \leq t \leq 2 \log 2$ 
with $L^1$-norm satisfying
\begin{equation} \label{ourODEsystem:d}
\|h\|_1\leq  c\,\left(\|\vec f\|_p+\|f_0\|_p\right).
\end{equation}
 Moreover, the new dependent variables $(\phi,\psi)$ are related to $(\vec v,\vec v_t)$ according to:
 \begin{equation} \label{ourODEsystem:e}
 \begin{pmatrix}
\phi \\ \psi
\end{pmatrix}
-\frac{1}{n^2}
 \begin{pmatrix}
n\vec v-\vec v_t \\ \vec v_t
\end{pmatrix}
\leq
c\,\e(t)\left(
|\vec v(t)|+|\vec v_t(t)|+\meanint|\nabla w|ds
\right).
\end{equation}
  \end{subequations}

Now, given $w$ with suitable properties, we solve (\ref{ourODEsystem:a}) with initial conditions $\phi(0)=0=\psi(0)$ to find $(\phi,\psi)$ and hence $\vec v$. 
But we want to separately control the dependence of $\vec v$ upon $w$, so let us write $\vec v=\vec v{\,}^w+\vec v{\,}^0$ where $\vec v{\,}^w$ corresponds to solving (\ref{ourODEsystem:a}) with $h(t)\equiv 0$ and $\vec v{\,}^0$ corresponds to solving (\ref{ourODEsystem:a}) with $g(t,\nabla w)\equiv 0$. We want to apply Proposition \ref{pr:Prop3} to estimate $\vec v{\,}^w$ on $(0,\infty)$, so we need to confirm that 
  $g=(g_1,g_2)$ satisfies $g_1\in L^1(0,\infty)$ and $g_2$ satisfies
(\ref{DS:g2}). 
To show that  $g_1\in L^1$, we use (\ref{ourODEsystem:c}) to conclude
 \[
 \begin{aligned}
 \int_0^\infty |g_1(t,\nabla w)|\,dt \leq & \ c\,\left(\int_0^\infty \e^2(t)\,dt\right)^{1/2}\left(\int_0^\infty\int_{S^{n-1}}|\nabla w|^2\,ds\,dt\right)^{1/2}\\
 & = c\,\left(\int_0^1\om^2(\rho)\frac{d\rho}{\rho}\right)^{1/2}\left(\int_0^1\meanint|\nabla w|^2\,ds\,\frac{d\rho}{\rho}\right)^{1/2},
 \end{aligned}
 \]
 and then invoke (\ref{eq:int-omega^2<delta}) to conclude
 \begin{equation}
\|g_1\|_1\leq c\, \sqrt\de\,\left(\int_0^1\meanint|\nabla w|^2\,ds\,\frac{d\rho}{\rho}\right)^{1/2}.
 \label{est:g_1}
\end{equation}
Similarly, to verify (\ref{DS:g2}), we estimate
 \begin{subequations}\label{est-g_2:joint}
\begin{equation}
 e^{\alpha t}\int_t^\infty |g_2(\tau)|\,e^{-\alpha \tau}\,d\tau
 \leq c\,\e(t)\int_t^\infty e^{\alpha(t-\tau)}\meanint |\nabla w|\,ds\,d\tau \leq c_\alpha \e(t),
   \label{est-g_2:a}
 \end{equation}
where 
\begin{equation}
c_\alpha=\frac {c}{\sqrt{2\alpha}}\,\left(\int_0^1\meanint |\nabla w|^2\,ds \frac{d\rho}{\rho}\right)^{1/2}.
 \label{est-g_2:b}
\end{equation}
 \end{subequations}
 As we shall see below, $M_2(\nabla w,r)=O(\om(r))$ as $r\to 0$, so the finiteness of $c_\alpha$ and the bound in (\ref{est:g_1}) follow from the following calculations:
for $j=0, 1, \dots$, let $r_j=2^{-j}$, so that
 \begin{subequations}\label{est-g_1:joint}
 \begin{equation}\label{est-g_1:a}
\int_0^1\int_{S^n}|\nabla w|^2\,ds\,\frac{d\rho}{\rho}
 =\sum_{j=0}^\infty \int_{r_j}^{2r_j}\int_{S^{n-1}} |\nabla w|^2\,ds\,\frac{d\rho}{\rho}
\leq c\,\sum_{j=0}^\infty M_2^2(\nabla w,r_j)  
\end{equation}
and
\begin{equation}\label{est-g_1:b}
\sum_{j=0}^\infty \om^2(r_j)=\sum_{j=0}^\infty \om^2(r_j)\frac{r_j-r_{j+1}}{r_{j+1}}
\leq c\int_0^1\om^2(\rho)\frac{d\rho}{\rho}\leq c\, \de,
 \end{equation}
  \end{subequations}
 where we have used (\ref{eq:int-omega^2<delta}) at the end. 
    
Now let us derive the PDE for $w$. Introduce the matrix $\Om=(\Om_{ij})$ with entries
  \begin{equation}
\Om_{ij}=a_{ij}-\de_{ij}, 
  \end{equation}
  and recall that $ |\Om_{ij}(r)|\leq \om(r)$ for $0<r<1$ and $ \Om_{ij}(r)=0$  for $r>1$.
  We can apply $I-P$ to (\ref{Lu=f:a}), to obtain
  \[
[\del_i(a_{ij}(x)\del_j u)]^\perp=[\del_i f_i+f_0]^\perp \qquad\hbox{in}\  \RR^n.
\]
If we  substitute (\ref{eq:spectral_decomp}) into this
 and  use $P[\Lap(u_0+\vec v\cdot \vec x)]=\Lap(u_0+\vec v\cdot \vec x)$ 
 and $ P[\Lap w]=0$, we  obtain the following equation for $w$:
 \[
 \Delta w+[ \hbox{div}(\Om\nabla w)]^\perp
+[ \hbox{div}(\Om\nabla (\vec v\cdot\vec x))]^\perp + [ \hbox{div}(\Om\nabla u_0)]^\perp
 =[\del_i f_i +f_0]^\perp.
  \]
  But we can use (\ref{ODE1-int:a}) to eliminate $u_0$ and  $\vec v=\vec v{\,}^w+\vec v{\,}^0$ to write this as: 
  \begin{equation}\label{eq:Lap-mu}
  \begin{aligned}
 \Delta w+&[ \hbox{div}(\Om\nabla w)]^\perp
 +[ \hbox{div}(\Om\nabla (\vec v{\,}^w\cdot\vec x))]^\perp
- [ \hbox{div}(\alpha^{-1}\Om\theta (r\vec\beta\cdot(\vec v{\,}^w)'+\vec \gamma\cdot\vec v{\,}^w+p[\nabla w] ) ) ]^\perp \\
 =&[\del_i f_i +f_0+\hbox{div}(\alpha^{-1}\Om\theta (r\vec\beta\cdot(\vec v{\,}^0)'+\vec \gamma\cdot\vec v{\,}^0)-\Om\nabla (\vec v{\,}^0\cdot\vec x))-\hbox{div}(\alpha^{-1}\Om\theta\vartheta)]^\perp.
 \end{aligned}
  \end{equation}
  
 Our strategy now is to simultaneously solve  (\ref{ourODEsystem:a}) and (\ref{eq:Lap-mu}) for $\vec v$ and $w$, then we can plug these into  (\ref{ODE1-int:a}) and integrate to find $u_0$ with $u_0(r)=0$ for $r>1/2$.  Indeed, for our chosen $p>n$, we assume $w$ is in the space $Y$ consisting of functions in $H^{1,p}_{\loc}(\RR^n\backslash\{0\})$ with finite norm
 \begin{equation}
 \|w\|_Y=\sup_{0<r<1}\frac{M_{1,p}(w,r)}{\om(r)\,r}+\sup_{r>1}\frac{M_{1,p}(w,r)}{\de\,r^{-n}}.
 \label{def:Y}
 \end{equation}
 We plug $w$ into (\ref{ourODEsystem:a}) and solve as described above to find $\vec v=\vec v{\,}^w+\vec v{\,}^0$; since
  $M_2(\nabla w,r)\leq c\,M_p(\nabla w,r)=O(\om(r))$ for $0<r<1$, we can apply Proposition 3 to estimate $\vec v{\,}^w$.
Now we want $w$ to satisfy (\ref{eq:Lap-mu}), so let us apply $\Lap^{-1}$, i.e.\ convolution by the fundamental solution,  to both sides of (\ref{eq:Lap-mu}) to obtain
 \begin{equation}
w+T[w]=\xi,
\label{eq:w+Tw}
\end{equation}
where
\[
 \begin{aligned}
 T[w]=  &\ \Lap^{-1}\left( [ \hbox{div}(\Om\nabla w)]^\perp
 +[ \hbox{div}(\Om\nabla (\vec v{\,}^w\cdot\vec x))]^\perp 
 - [ \hbox{div}(\alpha^{-1}\Om\theta (r\vec\beta\cdot(\vec v{\,}^w)'+\vec \gamma\cdot\vec v{\,}^w+p[\nabla w])]^\perp\right)\\
 \xi & =  \Lap^{-1}\left([\del_i f_i +f_0-\hbox{div}(\alpha^{-1}\Om\theta\vartheta)+\hbox{div}(\alpha^{-1}\Om\theta (r\vec\beta\cdot(\vec v{\,}^0)'+\vec \gamma\cdot\vec v{\,}^0)-\Om\nabla (\vec v{\,}^0\cdot\vec x))]^\perp\right).
 \end{aligned}
 \]
  We want to use the results of Sections 1 and 2 to show $\xi\in Y$ and $T$ is a bounded operator $Y\to Y$ with small norm.
  
To show $T:Y\to Y$ has small norm, consider $T[y]$ for $\|y\|_Y\leq 1$, i.e.\ we assume $y$ satisfies
\begin{subequations}\label{M1p(w):joint}
\begin{equation}
M_{1,p}(y,r)\leq \om(r)\,r
\quad\hbox{for}\ 0<r<1,
\label{M1p(w):0<r<1}
\end{equation}
and
\begin{equation}
M_{1,p}(y,r)\leq \de \,r^{-n}=\om(r)\,r^{-n}
\quad\hbox{for}\ r>1.
\label{M1p(w):r>1}
\end{equation}
\end{subequations}
Now let us consider separately the three terms,
\[
\begin{aligned}
 T_{1}[y]= &  \Lap^{-1}\left( [ \hbox{div}(\Om\nabla y)]^\perp\right), \quad
 T_{2}[y]= \Lap^{-1}\left([ \hbox{div}(\Om\nabla (\vec v{\,}^y\cdot\vec x))]^\perp\right),\\
T_{3}[y] & =   \ \Lap^{-1}\left([ \hbox{div}(\alpha^{-1}\Om\theta (r\vec\beta\cdot(\vec v{\,}^y)'+\vec \gamma\cdot\vec v{\,}^y+p[\nabla y]))]^\perp\right).
\end{aligned}
 \]

First we consider $T_1$. 
Using Proposition 1 and (\ref{eq:Mp(omega_f)}), we have
\begin{equation}
M_{1,p}(T_1[y],r)\leq c\,\left(r^{-n}\int_0^r M_p(\Om\nabla y,\rho)\rho^n\,d\rho
+r^2\int_r^\infty M_p(\Om\nabla y,\rho)\rho^{-2}\,d\rho\right).
\label{est:M1p(T1).1}
\end{equation}
But recall $|\Om(r)|\leq \om(r)$ for $0<r<1$ and $\Om(r)=0$ for $r>1$. Thus, for $0<r<1$,  
we can use  (\ref{M1p(w):0<r<1}) and  (\ref{M1p(w):r>1}) 
with the facts that $\om(r)$ is nondecreasing and $\om(r)r^{-1+\kappa}$ is nonincreasing
to obtain
\begin{subequations}\label{M1p(T1):joint}
\begin{equation}
\begin{aligned}
M_{1,p}(T_1[y],r) & \leq  c\,\left(r^{-n}\int_0^r\om^2(\rho)\rho^n\,d\rho
+r^2\int_r^1\om^2(\rho)\rho^{-2}\,d\rho\right)\\
\leq \ & c\,\left(r^{-n}\om^2(r)r^{n+1}+r^2\de\om(r)r^{-1+\kappa}(r^{-\kappa}-1)\right)\leq c\,\de\,\om(r)\,r.
\end{aligned}
\label{M1p(T1):a}
\end{equation}
For $r>1$ we simply get
\begin{equation}
M_{1,p}(T_1[y],r)\leq c\,r^{-n}\int_0^1\de^2\rho^n\,d\rho\leq c\,\de^2\,r^{-n}.
\label{M1p(T1):b}
\end{equation}
\end{subequations}
The estimates   (\ref{M1p(T1):a}) and (\ref{M1p(T1):b}) show that $\|T_1[y]\|_Y\leq c\,\de$, so
for $\de$ sufficiently small we conclude that $T_1:Y\to Y$ has norm less than $1/3$.

Next we consider $T_2$. To use Proposition 1, we need to estimate $M_p(\nabla( \vec v{\,}^y\cdot\vec x),r)$. 
But $\vec v{\,}^y$ and $r(\vec v{\,}^y)'$ can be expressed in terms of the solutions $\phi,\psi$ of (\ref{ourODEsystem:a}) with $h\equiv 0$, to which we can apply Proposition 3 to find 
\begin{equation}\label{est:grad-v}
\sup_{|x|<1}|\nabla(\vec v{\,}^y\cdot \vec x)|\leq 
c\,\sup_{r<1} \,(r|(\vec v{\,}^y)'| + |\vec v{\,}^y|)\leq 
c\,\sup_{t>0}\,(|\phi|+|\psi|)\leq
c\,(c_\alpha+\|g_1\|_1),
\end{equation}
where $c_\alpha$ is given in (\ref{est-g_2:b}) and $\|g_1\|$ is estimated as in (\ref{est:g_1}). But, using (\ref{est-g_1:a}), $p>2$, (\ref{M1p(w):0<r<1}), and (\ref{est-g_1:b}), we have
\[
\int_0^1\meanint_{S^{n-1}}|\nabla y|^2\,ds\,\frac{d\rho}{\rho}\leq c\,\sum_{j=0}^\infty M_2^2(\nabla y,r_j)
\leq c\,\sum_{j=0}^\infty M_p^2(\nabla y,r_j)\leq c\,\sum_{j=0}^\infty \om^2(r_j)\leq c\,\de.
\]
We conclude that for $0<r<1$:
\begin{equation}\label{est:Mp(nabla_vw)}
M_p(\nabla( \vec v{\,}^y\cdot \vec x),r)\leq 
c\,\sup_{|x|<1}|\nabla(\vec v{\,}^y\cdot\vec x)|\leq
c\,\sqrt\de.
\end{equation}
We can use this in Proposition 1 and estimate as above to obtain for $0<r<1$:
\begin{subequations}\label{M1p(T3):joint}
 \begin{equation}
 M_{1,p}(T_2[y],r)  \, 
 \leq \ c\left(r^{-n}\int_0^r\om(\rho)\sqrt\de\rho^n\,d\rho
 +r^2\int_r^1\om(\rho)\sqrt\de\rho^{-2}\,d\rho\right)
 \leq c\,\om(r)\,\sqrt\de\,r.
 \label{M1p(T3):a}
\end{equation}
 Meanwhile, for $r>1$ we know $\vec v=\vec v{\,}^w\equiv 0$, so by Proposition 1
 \begin{equation}
M_{1,p}(T_2[y],r)\leq c\,r^{-n}\int_0^1 \om(\rho)\,M_p(\nabla(\vec v{\,}^y\cdot\vec x),\rho)\rho^n\,d\rho \leq
c\,\de^{3/2}\,r^{-n}.
\label{M1p(T3):b}
\end{equation}
\end{subequations}
For $\de$ sufficiently small, we conclude from (\ref{M1p(T3):a}), (\ref{M1p(T3):b}) that  $T_2:Y\to Y$ has norm less than $1/3$.

To show $T_3:Y\to Y$ is small, we need to estimate $M_p(\alpha^{-1}\Om\theta (r\vec\beta\cdot(\vec v{\,}^y)'+\vec \gamma\cdot\vec v{\,}^y+p[\nabla y]),r)$ only for $0<r<1$ (since it vanishes for $r>1$). But recalling the bounds on $\Omega$, $\vec \beta$, $\vec\gamma$ and 
$p[\nabla y]$, we have for $0<r<1$:
\[
\begin{aligned}
M_p(\alpha^{-1}\Om\theta (r\vec\beta\cdot(\vec v{\,}^y)'+\vec \gamma\cdot \vec v{\,}^y)+p[\nabla y],r)
&\leq c\,\om^2(r)\left(M_p(r(\vec v{\,}^y)',r)+M_p(\vec v{\,}^y,r)+M_p(\nabla y,r)
\right)\\
&\leq c\,\sqrt{\de}\,\om^2(r),
\end{aligned}
\]
where at the end we have used Proposition 3 and (\ref{est:grad-v}), similar to our derivation of (\ref{est:Mp(nabla_vw)}).
Applying Proposition 1,  we obtain for $0<r<1$
\[
M_{1,p}(T_3[y],r)
\leq c\left(r^{-n}\int_0^r \sqrt\de\,\om^2(\rho)\,\rho^n\,d\rho+r^2\int_r^1\sqrt\de\,\om^2(\rho)\,\rho^{-2}\,d\rho\right) \leq c\,\de^{3/2}\,\om(r)\,r,
\]
and for $r>1$ we have simply
\[
\begin{aligned}
M_{1,p}(T_3[y],r)& \leq c\,r^{-n}\int_0^1\sqrt{\de}\,\om(\rho)\rho^n\,d\rho \leq c\,\de^{3/2}\, r^{-n}.
\end{aligned}
\]
Taking $\de$ sufficiently small, we conclude that $T_3:Y\to Y$ has norm less than 1/3.
We have therefore shown that $T=T_1+T_2+T_3:Y\to Y$ has norm less than $1$.

  To show $\xi\in Y$, first note that ${\rm supp}\ f\subset A_{1/4}=\{x:1/4\leq |x|\leq 1/2\}$ implies $M_p(f,r)=0$ for $r<1/8$ and $r>1/2$, whereas  $M_p(f,r)\leq c\,\|f\|_p$ for $1/8\leq r\leq 1/2$. Now we separately treat
  \[
  \xi_1=\Delta^{-1}([\del_i f_i]^\perp), \quad
  \xi_2=\Delta^{-1}([f_0]^\perp), \quad
  \xi_3=\Delta^{-1}([\hbox{div}(\alpha^{-1}\Om\theta\vartheta(r))]^\perp),
  \]
  and
  \[
  \xi_4=\Delta^{-1}([\hbox{div}(\alpha^{-1}\Om\theta (r\vec\beta\cdot(\vec v{\,}^0)'+\vec \gamma\cdot\vec v{\,}^0)-\Om\nabla (\vec v{\,}^0\cdot\vec x))]^\perp).
  \]
Since both $f_i$ and $f_0$ are supported in $A_{1/4}$, we can apply Proposition 1 to estimate $\xi_1$ as
    \[
M_{1,p}(\xi_1,r)   \leq   \begin{cases}
 c\,\|\vec f\|_p\,r^2 & \hbox{for}\ 0<r<1,\\
 c\,\|\vec f\|_p\,r^{-n} &  \hbox{for}\  r>1,
\end{cases}
\]
and Proposition 2 to estimate $\xi_2$ as
    \[
M_{1,p}(\xi_2,r)\leq M_{2,p}(\xi_2,r)   \leq   \begin{cases}
 c\,\| f_0\|_p\,r^2 & \hbox{for}\ 0<r<1,\\
 c\,\| f_0\|_p\,r^{-n} &  \hbox{for}\  r>1.
\end{cases}
\]
Since $\Om\theta\vartheta=0$ for $|x|=r>1$ and $r<1/8$, whereas $M_p(\alpha^{-1}\Omega\theta\vartheta,r)\leq c\,(\|\vec f\|_p+\|f_0\|_p)$ for $1/8<r<1$, we similarly conclude
\[
   M_{1,p}(\xi_3,r)\leq    \begin{cases}
 c\,(\|\vec f\|_p+\| f_0\|_p)\,r^2 & \hbox{for}\ 0<r<1,\\
 c\,(\|\vec f\|_p+\| f_0\|_p)\,r^{-n} &  \hbox{for}\  r>1.
\end{cases}
\]
To estimate $\xi_4$ we need to estimate $M_p(\nabla( \vec v{\,}^0\cdot\vec x),r)$. 
But $\vec v{\,}^0$ and $r(\vec v{\,}^0)'$ can be expressed in terms of the solutions $\phi,\psi$ of (\ref{ourODEsystem:a}) with $g\equiv 0$, to which we can apply Proposition 3 to find 
\[
\sup_{|x|<1}|\nabla(\vec v{\,}^0\cdot \vec x)|\leq c\,\|h\|_1.
\]
Combined with (\ref{ourODEsystem:d}), we find $M_p(\nabla( \vec v{\,}^0\cdot\vec x),r)\leq c\,(\|\vec f\|_p+\|f_0\|_p)$,
and then applying Proposition 1 yields
\[
\begin{aligned}
  M_{1,p}(\xi_4,r)\leq &
\,  c\,(\|\vec f\|_p+\|f_0\|_p)\left(r^{-n}\int_0^r\om(\rho)\rho^n\,d\rho+r^2\int_r^1\om(\rho)\,\rho^{-2}\,d\rho\right)\\
 \leq &  \begin{cases}
 c\,(\|\vec f\|_p+\| f_0\|_p)\,\om(r)\,r & \hbox{for}\ 0<r<1,\\
 c\,(\|\vec f\|_p+\| f_0\|_p)\,r^{-n} &  \hbox{for}\  r>1.
\end{cases}
  \end{aligned}
\]
Plugging these estimates  into (\ref{def:Y}), we find that $\xi\in Y$ with 
\[
\|\xi\|_Y\leq c\, (\|\vec f\|_p+\|f_0\|_p).
\]
 
We conclude that (\ref{eq:w+Tw}) admits a unique solution $w\in Y$ satisfying 
\begin{equation}\label{est:w_Y}
\|w\|_Y\leq c\, (\|\vec f\|_p+\|f_0\|_p),
\end{equation}
and then we use this to find $\vec v^w$ and $u_0^w$ as described above. We also know from Proposition 1 that $Pw=0$, which in particular shows that $\int_{|x|<r}w(x)\,dx=0$ for every $r>0$. Since $p>n$, we can apply Morrey's inequality  (cf.\ Theorem 7.17 in \cite{GT}) to obtain
\begin{equation}
\sup_{|x|<r}|w(x)|\leq c_n\,r\,\left(\meanint_{|y|<r}|\nabla w|^p\,dy\right)^{1/p}.
\label{MorreyIneq}
\end{equation}
But for fixed $r\in (0,1)$,  we can introduce $r_j=r\,2^{-j}$  to compute
\[
\begin{aligned}
\meanint_{|y|<r}|\nabla w|^p\,dy=\sum_{j=0}^\infty 2^{-jn}r_j^{-n}\int_{r_{j+1}<|y|<r_j}|\nabla w|^p\,dy\leq c\sup_{0<\rho<1}M_p(\nabla w,\rho).
\end{aligned}
\]
Recalling that (\ref{est:w_Y}) implies $M_p(\nabla w,r)\leq c\,\om(r)\,(\|\vec f\|_p+\|f_0\|_p)$ for $0<r<1$, we find
\begin{equation}\label{est:sup-w}
\sup_{|x|<r}|w(x)|\leq c\,r\,\om(r)\,(\|\vec f\|_p+\|f_0\|_p)\quad\hbox{for}\ 0<r<1.
\end{equation}
In particular, this implies that $w$ is differentiable at $x=0$ with $\del_jw(0)=0$. 

What about $\vec v$ and $u_0$? Since we now know $M_p(\nabla w,r)\leq c\,\om(r)\,(\|\vec f\|_p+\|f_0\|_p)$, we obtain
\[
\int_0^1\meanint |\nabla w|^2\,ds\,\frac{d\rho}{\rho}\leq c\,\sum_{j=0}^\infty M_p^2(\nabla w,r_j)
\leq c\,(\|\vec f\|_p+\|f_0\|_p)^2\sum_{j=0}^\infty \om^2(r_j)\leq c\,\de\,(\|\vec f\|_p+\|f_0\|_p)^2.
\]
Consequently, our analysis of (\ref{ourODEsystem:a}) yields
\begin{equation}\label{est:sup-v}
\sup_{0<r<1}(|\vec v(r)|+r\,|\vec v{\,}'(r)|)\leq c\,(\|\vec f\|_p+\|f_0\|_p).
\end{equation}
Using (\ref{ODE1:a}), we perform the following estimates
\[
|u_0(r)-u_0(0)|  \leq \int_0^r |u_0'(\rho)|\,d\rho
 \leq c\,\int_0^r \left( |\vartheta(\rho)|+|\rho\vec\beta(\rho)\cdot\vec v{\,}'(\rho)|+|\vec\gamma(\rho)\cdot\vec v|+\om(\rho)\meanint|\nabla w|\,ds\right)\,d\rho
\]
\[
\int_0^r \left( |\vartheta(\rho)|+|\rho\vec\beta(\rho)\cdot\vec v{\,}'(\rho)|+|\vec\gamma(\rho)\cdot\vec v|\right)d\rho
\leq c\,\om(r)\,r\,\left(\|\vec f\|_p+\|f_0\|\right)
\]
\[
\int_0^ r\om(\rho)\meanint|\nabla w|\,ds\,d\rho\leq \om(r)\,r\int_0^1\meanint|\nabla w|\,ds\,\frac{d\rho}{\rho}
\leq c\,\om(r)\,r\,\left(\|\vec f\|_p+\|f_0\|\right)
\]
 to obtain
\begin{equation}
|u_0(r)-u_0(0)|\leq c\,\om(r)\,r\,(\|\vec f\|_p+\|f_0\|_p).
\label{eq:u0(r)-u0(0)}
\end{equation}
But (\ref{eq:u0(r)-u0(0)}) implies that $u_0$ is differentiable at $r=0$ with $u_0'(0)=0$.

Thus we have found a solution of (\ref{Lu=f:a}) in the form $\tilde u(x)=u_0(r)+\vec v(r)\cdot\vec x+w(x)$ which is Lipschitz continuous at $x=0$, but we need to verify that $\tilde u$ coincides with the solution $u$ of (\ref{Lu=f:a}) that we began with. However, if we let $z(x)=u(x)-\tilde u(x)$, we find that $z$ is a weak solution of the homogenous equation ${\mathcal L}z=0$ in $\RR^n$, and $z(x)\to 0$ as $x\to \infty$ (since $|x|>1$ implies $u(x)\equiv 0$ and $\tilde u(x)=w(x)=O(|x|^{-n})$ as $|x|\to\infty$). The maximum principle shows that $z\equiv 0$, i.e.\ $u=\tilde u$.

We conclude that our solution $u$ of (\ref{eq:Lu=0}) is Lipschitz continuous at $x=0$. To obtain the desired estimate (\ref{est:LipBnd}), we first
combine (\ref{est:sup-w}), (\ref{est:sup-v}),  (\ref{eq:u0(r)-u0(0)}), and recall the definitions of $\vec f$ and $f_0$ to conclude 
\begin{equation}\label{est:u(x)-u(0)}
|u(x)-u(0)|\leq c\,|x|\,(\|\nabla u\|_{L^p(B')}+\|u\|_{L^p(B')})\quad\hbox{for}\ 0<|x|<1/2\quad\hbox{and}\ B'=B_{1/2}(0).
\end{equation} 
But as a solution of (\ref{eq:Lu=0}), $u$ satisfies the elliptic estimate 
\begin{equation}\label{est:ADN}
\|\nabla u\|_{L^p(B')}\leq c\,\|u\|_{L^p(B^*)}, \quad\hbox{where} \ B^*=B_{3/4}(0),
\end{equation}
which  can be found, for example, in \cite{ADN} when the coefficients are continuous; however, their proof extends directly to the case where the coefficients have small oscillation, which we may assume in the unit ball by taking $\de$ sufficiently small. But from \cite{Mo}, $u$ also satisfies the following estimate:
\begin{equation}\label{est:Moser}
\sup_{|y|\leq 3/4}|u(y)|\leq c\,\|u\|_{L^2(B)}, \quad\hbox{where} \ B=B_1(0).
\end{equation}
Using these in (\ref{est:u(x)-u(0)}), we obtain 
\begin{equation}\label{est-improved:u(x)-u(0)}
|u(x)-u(0)|\leq  c\,|x|\|u\|_{L^2(B)}
\quad\hbox{for}\ 0<|x|<1/2,
\end{equation} 
which is  (\ref{est:LipBnd}) for $r=1$. The case of general $r\in (0,1)$ can be achieved by scaling: $\tilde x=x/r$ and $\tilde u(\tilde x)=u(x)$. Thus (\ref{est:LipBnd}) is proved.

Now let us add the hypothesis that every solution of (\ref{eq:DynSystm}) is asymptotically constant. Then, according to Proposition 3, the solution $(\phi,\psi)$ satisfies $\phi(t)\to \phi_\infty$ and $\psi(t)\to 0$ as $t\to \infty$. Using (\ref{ourODEsystem:e}), we see that $\vec v_t\to 0$ as $t\to \infty$ and consequently
\[
\lim_{r\to 0}\vec v(r)=n\lim_{t\to\infty}\phi(t)=n\,\phi_\infty
\quad\hbox{and}\quad
\lim_{r\to 0}r\vec v{\,}'(r)=0.
\]
 In particular, $\vec v(|x|)\cdot \vec x$ is differentiable at $\vec x=0$.
Putting this together with what we have already found about $u_0$ and $w$, we conclude that our weak solution $u$ of (\ref{eq:Lu=0}) is differentiable at $x=0$ and (\ref{eq:del_j_u(0)}) holds.
\hfill$\Box$

 \section{Proof of Corollaries}
 
To prove Corollary 1, let us write (\ref{eq:DynSystm}) as
\begin{equation}
\dot \phi = B(t)\phi\quad\hbox{for}\ T<t<\infty,
\label{eq:ODE}
\end{equation}
where $B(t)$ is bounded, but not necessarily self-adjoint. (In fact, we know $B(t)\to 0$ as $t\to\infty$, but we will not need this fact here.) Let $\mu(B)$ denote the largest eigenvalue of $(B+B^t)/2$, which satisfies (cf.\ \cite{C}, Ch.II, Sec.1)
\[
\mu(B)=\lim_{h\to 0^+}\frac{\pmb{\bigr |}I+hB\pmb{\bigr |}-1}{h}.
\]
Assuming that $B(t)$ is continuous, the following inequality is proved in \cite{C} (Ch.III, Sec.2):
\begin{equation}
|\phi(t)|\leq |\phi(s)|\exp\left(\int_s^t\mu(B(\tau))\,d\tau\right)
\quad\hbox{for}\ t\geq s\geq T.
\label{eq:unif-stable-cond}
\end{equation}
Let us verify that (\ref{eq:unif-stable-cond}) holds even though $B(t)$ may be discontinuous.

Let $f(t)=|\phi(t)|$ for $t\geq T$. It is easy to see from (\ref{eq:ODE}) that $\phi$ is Lipschitz continuous, so Rademacher's theorem implies that $\phi$ is differentiable almost everywhere. Consequently, almost everywhere $f$ has a right-hand derivative $\dot f_+$ satisfying
\[
\dot f_+(t)=\lim_{h\to 0^+} \frac{|\phi(t)+h\dot\phi(t)|-|\phi(t)|}{h}=\lim_{h\to 0^+} \frac{|\phi(t)+hB(t)\phi(t)|-|\phi(t)|}{h}\leq \mu(B(t))\,f(t)
\quad\hbox{a.e}.
\]
Since $\mu(B(t))$ is bounded, we know that 
\[
w(t):=f(t)\,\exp\left(-\int_{T}^t\mu(B(\tau))\,d\tau\right)
\]
is continuous and has a right-hand derivative satisfying
\[
\dot w_+(t)=\left[\dot f_+(t)-\mu(B(t))f(t)\right]\exp\left(-\int_{T}^t\mu(B(\tau))\,d\tau\right)\leq 0\quad\hbox{a.e}.
\]
We conclude that $w(t)$ is nonincreasing for almost every $t>T$. But $w(t)$ is continuous, so $w(t)$ is nonincreasing for all $t>T$. We conclude that  (\ref{eq:unif-stable-cond}) holds. Moreover, (\ref{eq:unif-stable-cond}) together with
\[
\int_s^t \mu(B(\tau))\,d\tau <K \quad\hbox{for} \ t>s>T
\]
implies that (\ref{eq:ODE}) is uniformly stable. Thus we may apply Theorem 1 to obtain Corollary 1.

To prove Corollary 2, we assume (\ref{Differentiability:joint}) and introduce a change of dependent variable  (as in \cite{E} Section 11.1):
\[
\phi(t)=(I+ S(t))\,\xi(t)
\]
where
\[
 S(t)=\int_t^\infty  \widetilde R(\tau)\,d\tau, \quad\hbox{with}\ \widetilde R(t)=   R(e^{-t}).
\]
We find that $\xi$ satisfies
\[
(I+ S(t))\frac{d\xi}{dt}+ \widetilde R S\,\xi=0.
\]
But $ S(t)\to 0$ as $t\to\infty$, so we can take $T$ sufficiently large and conclude that $(I+S(t))$ is invertible, and consequently $\xi$ satisfies
\begin{equation}
\frac{d\xi}{dt}+ Q(t)\,\xi=0,
\label{eq:newDynSystm}
\end{equation}
where by hypothesis we have $ Q=(1+ S)^{-1} \widetilde R S\in L^1(T,\infty)$. As we have observed, 
$ Q\in L^1(T,\infty)$ implies that all solutions
of (\ref{eq:newDynSystm}) are asymptotically constant, i.e.\ $\xi(t)= \xi_\infty+o(1)$ as $t\to\infty$, so
\[
\phi(t)=(I+S(t))(\xi_\infty+o(1))=\xi_\infty+o(1)
\]
and we see that all solutions of (\ref{eq:DynSystm}) are asymptotically constant. We have also observed that $ Q\in L^1(T,\infty)$ implies that (\ref{eq:newDynSystm}) is uniformly stable, so the same is true of  (\ref{eq:DynSystm}). 

To prove Corollary 3, observe that $\int_T^t \mu(B(\tau))\,d\tau\to -\infty$ as $t\to \infty$ together with (\ref{eq:unif-stable-cond}) implies that all solutions of (\ref{eq:ODE}) tend to zero as $t\to\infty$, i.e.\ the null solution is asymptotically stable. But this implies that $\vec v(r)$ in (\ref{eq:Taylor1}) satisfies $\vec v(r)\to 0$ as $r\to 0$ and Corollary 3 follows.

\section{Examples of Gilbarg-Serrin Type}

In \cite{GS}, Gilbarg and Serrin consider examples of the form
\begin{equation}\label{eq:GSexample}
a_{ij}(x)=\de_{ij}+g(r)\theta_i\theta_j, 
 \end{equation}
 where $g(0)=0$ but vanishes slowly as $r\to 0$. They use such examples to show that Dini continuity is essential for their ``extended maximum principle'' to hold, but we shall use them to explore the conditions in Theorem 1 and its corollaries. 
We assume that $|g(r)|\leq \om(r)$ for $r$ near $0$ with $\om$ satisfying (\ref{eq:Sq-Dini}), and we can explicitly calculate the quantities introduced in the Appendix:
\[
\alpha(r)=1+g(r),\quad \beta(r)=0=\gamma(r), \quad A(r)=B(r)=\frac{1+g(r)}{n}\,I, \quad 
C(r)=\left(1+\frac{g(r)}{n}\right)I.
 \]
Moreover, the matrix (\ref{def:R}) is given by
\[
R_{\ell k}(r)=\frac{1-n}{n}\, g(r)\, \de_{\ell k},
\]
so the dynamical system (\ref{eq:DynSystm}) reduces to the scalar equation
\begin{equation}
\frac{d\phi}{dt}=\frac{n-1}{n}\,  \tilde g(t)\,\phi,
\label{eq:g-ODE}
\end{equation}
where $\tilde g(t)=g(e^{-t})$. 

Now consider a weak solution $u$ of (\ref{eq:Lu=0}) in a domain containing $x=0$ and $a_{ij}$ of the form (\ref{eq:GSexample}). According to Theorem 1, $u$ is Lipschitz continuous at $x=0$ provided (\ref{eq:g-ODE}) is uniformly stable for $t>T$ with $T$ sufficiently large; but
it is easy to solve (\ref{eq:g-ODE}) and see that it is uniformly stable if and only if
\begin{equation}
\int_s^t  \tilde g(\tau)\,d\tau < K \quad\hbox{for}\ t>s>T.
\label{eq:g-unifstble}
\end{equation}
Moreover, $\mu(\mathcal{R}(r))=(1-n^{-1})g(r)$, so (\ref{eq:Lipschitz-condition}) agrees with (\ref{eq:g-unifstble}) and we see that Corollary 1 is sharp for this class of examples. On the other hand, solutions of (\ref{eq:g-ODE}) are asymptotically constant if and only if the improper integral
\begin{equation}
\int_T^\infty  \tilde g(\tau)\,d\tau \hbox{ converges to an extended real number} < \infty.
\label{eq:g-asym-const}
\end{equation}
Thus  Theorem 1 implies that $u$ is differentiable at $x=0$ if both (\ref{eq:g-unifstble}) and (\ref{eq:g-asym-const}) hold. The case
\begin{equation}
\int_T^\infty \tilde g(t)\,dt=-\infty
\label{eq:int-g=-infty}
\end{equation}
 in (\ref{eq:g-asym-const}) pertains to Corollary 3, which is sharp for this class of examples. On the other hand, the case that $\tilde g(t)$ is integrable pertains to Corollary 2 and coincides  with the hypothesis (\ref{Differentiability:1}); however, in Corollary 2 we also require (\ref{Differentiability:2}), since (as mentioned in the Introduction) the condition (\ref{Differentiability:1}) alone does not imply the uniform stability of (\ref{eq:DynSystm}). 

In fact, this class of examples may be used to show not only this last statement, but in general  that uniform stability is not implied by every solution being asymptotically constant: we only need to construct $\tilde g(t)$ for which (\ref{eq:g-asym-const}) holds but (\ref{eq:g-unifstble}) fails. Moreover, if  (\ref{eq:g-asym-const}) holds because $\tilde g$ is integrable on $(0,\infty)$, then this explains the need for a condition such as (\ref{Differentiability:2}) in Corollary 2; on the other hand, if we construct $\tilde g$ for which
(\ref{eq:int-g=-infty}) holds 
 and yet (\ref{eq:g-unifstble}) fails, then we see that
 (\ref{eq:Lipschitz-condition}) is not implied by
(\ref{int-Rgoesto-infty}), so both conditions are necessary in Corollary 3. In this latter regard, let us observe that
 \cite{C} gives an explicit example of a function $\tilde g(t)$ satisfying (\ref{eq:int-g=-infty}) and yet (\ref{eq:g-unifstble}) fails: there exist $t_j\to\infty$ for which $\int_{t_{2j}}^{t_{2j+1}}  \tilde g(\tau)\,d\tau\to\infty$ and yet $\int_{t_{2j+1}}^{t_{2j+2}} \tilde g(\tau)\,d\tau\to -\infty$ more rapidly so that (\ref{eq:int-g=-infty}) still holds. Now the example in \cite{C} does not have $\tilde g(t)\to 0$ as $t\to \infty$, but the example can be modified to achieve this; in fact, we can even arrange $\tilde g(t)=O(t^{-2/3})$, which implies that $\tilde g\in L^2(T,\infty)$ and so the $a_{ij}$ are square-Dini continuous at $x=0$. Moreover, the example can be modified so that (\ref{eq:int-g=-infty}) is replaced by the condition that $\tilde g$ is integrable on $(0,\infty)$. Thus (\ref{eq:g-unifstble}) and (\ref{eq:g-asym-const}) are completely independent conditions, even under the assumption that the  coefficients $a_{ij}$ are square-Dini continuous at $x=0$.
 
\section{Appendix}
In this appendix we provide the details behind the derivation of the dynamical system (\ref{ourODEsystem:joint}).
To express this system, let us introduce
$n\times n$ matrices $A, B, C$ and vectors $\vec \xi$, $\vec \zeta$ by
\begin{align}\notag
A_{\ell k}(r)=\meanint_{S^{n-1}} a_{ij}(r\theta)\theta_i\theta_j\theta_\ell\theta_kds_\theta,\ 
B_{\ell k}(r)&=\meanint_{S^{n-1}} a_{\ell j}(r\theta)\theta_j\theta_kds_\theta,\ 
C_{\ell j}(r)=\meanint_{S^{n-1}}  a_{\ell j}(r\theta)\,ds_\theta,\\ \notag
\xi_\ell[\nabla w](r)=\meanint_{S^{n-1}} a_{ij}\theta_i\theta_\ell\del_jw\,ds_\theta,
\quad
&\zeta_\ell[\nabla w](r)=\meanint_{S^{n-1}} a_{\ell j}\del_jw \,ds_\theta,\notag
\end{align}
which satisfy for $0<r<1$
\begin{equation}
\begin{aligned}
\pmb{|}\,A-n^{-1}I\,\pmb{|},\ \pmb{|}\,B-n^{-1}I\,\pmb{|},\ \pmb{|}\,C-I\,\pmb{|} \leq \om(r)\\
|\vec\xi[\nabla w](r)|,\  |\vec\zeta[\nabla w](r)|
\leq \om(r)\,\meanint_{S^{n-1}} |\nabla w|\,ds,
\end{aligned}
\label{est:A-etc}
\end{equation}
while for $r>1$ we use (\ref{eq:aij_for_r>1}) to conclude
$A=n^{-1}I=B$, $C=I. $ (Notice that the matrix $A(r)$ introduced above is not the same as the matrix $A(x)$ used in the Introduction.)
Now using $\eta=\eta(r)x_\ell$ in (\ref{def:Lu=f}), we obtain 
the 2nd-order system of ODEs
\begin{equation}
-\left[ r^n(u_0'\vec\beta+rA\vec v{\,}'+B\vec v+\vec\xi[\nabla w]-\vec f{\,}^\#)\right]'+
r^{n-1}(u_0'\vec\gamma+rB\vec v{\,}'+C\vec v+\vec\zeta[\nabla w]+\vec f{\,}^\flat)=0,
\label{eq:ODE2}
\end{equation}
where 
\[
f^\#_{\ell}(r)=\meanint_{S^{n-1}}f_i(r\theta)\theta_i\theta_\ell\,ds_\theta \quad\hbox{and}\quad
f^\flat_{\ell}(r)=\meanint_{S^{n-1}}f_0(r\theta)\theta_\ell\,ds_\theta
\]
are supported in $1/4\leq r\leq 1/2$.
Next we can use (\ref{ODE1:a}) to solve for $u_0'(r)$ and eliminate $u_0$ from (\ref{eq:ODE2}); this is  the 2nd-order system of ODEs (depending upon $w$) that we want to analyze:
\begin{equation}\label{eq:ODE3}
\begin{aligned}
-&r\left[r^{n}\left(rA\vec v{\,}' + B\vec v+\vec\xi[\nabla w]-\vec f{\,}^\#\right)-\frac{r^{n}}{\alpha(r)}\left(r\vec\beta\cdot\vec v{\,}'+\vec\gamma\cdot\vec v+p[\nabla w]-\vartheta\right)\vec\beta\right]' \\ 
&+
\left[r^{n}\left(rB\vec v{\,}'+C\vec v+\vec\zeta[\nabla w]+\vec f{\,}^\flat\right)-\frac{r^{n}}{\alpha(r)}\left(r\vec\beta\cdot\vec v{\,}'+\vec\gamma\cdot\vec v+p[\nabla w]-\vartheta\right)\vec\gamma\right]=0.
\end{aligned}
\end{equation}
If we make the substitution $r=e^{-t}$, then (\ref{eq:ODE3}) becomes 
\[
\begin{aligned}
&\left[e^{-nt}\left(-A\vec v_t+B\vec v+\vec\xi[\nabla w]-\vec f{\,}^\#-\frac{1}{\alpha}(-\vec\beta\cdot\vec v_t+\vec\gamma\cdot\vec v+p[\nabla w]-\vartheta)\vec\beta\right)\right]_t \cr
+&\ e^{-nt}\left(-B\vec v_t+C\vec v+\vec\zeta[\nabla w]+\vec f{\,}^\flat-\frac{1}{\alpha}(-\vec\beta\cdot\vec v_t+\vec\gamma\cdot\vec v+p[\nabla w]-\vartheta)\vec\gamma\right)=0,
\end{aligned}
\]
which after some rearrangement can be written
\[
\begin{aligned}
\left[-A\vec v_t+B\vec v+\vec\xi[\nabla w]-\vec f{\,}^\#+\frac{\vec\beta\cdot\vec v_t-\vec\gamma\cdot\vec v-p[\nabla w]+\vartheta}{\alpha}\vec\beta\right]_t -(B-nA)\vec v_t+\frac{\vec\beta\cdot\vec v_t}{\alpha}(\vec\gamma-n\vec\beta)&
 \cr
+(C-nB)\vec v-\frac{\vec\gamma\cdot\vec v}{\alpha}(\vec\gamma-n\vec\beta)
=n\left[\vec\xi[\nabla w]-\frac{p[\nabla w]-\vartheta}{\alpha}\vec\beta\right]
+\frac{p[\nabla w]-\vartheta}{\alpha}\vec\gamma
-\vec\zeta[\nabla w]-\vec f{\,}^\flat&
\end{aligned}
\]

To avoid differentiating the coefficient matrices, let us convert this to a 1st-order system for the $2n$-vector $V=(V_1,V_2)$ where $V_1=\vec v$ and
\begin{equation}\label{def:V2}
V_2:=-A\vec v_t+B\vec v+\vec \xi[\nabla w] -\vec f{\,}^\#+\frac{\vec\beta\cdot\vec v_t-\vec\gamma\cdot\vec v-p[\nabla w]+\vartheta}{\alpha}\vec\beta.
\end{equation}
Notice that, as in Section 2, we have omitted arrows over the vectors $V_1$, $V_2$, and $V$, and we shall now also use the dot notation for $d/dt$. Since $A$ is invertible (for $\de$ sufficiently small), our 1st-order system can be written as
\begin{equation}\label{eq:ODE6}
\begin{aligned}
\dot V_1-&A^{-1}B V_1+A^{-1} V_2 + \left[\frac{\vec\beta\cdot\dot V_1-\vec\gamma\cdot V_1}{\alpha}\right]A^{-1} \vec\beta =A^{-1}\left[\vec\xi[\nabla w]-\vec f{\,}^\#-\frac{p[\nabla w]-\vartheta}{\alpha}\vec\beta\right]\cr
\dot V_2+&(C-BA^{-1}B) V_1+(BA^{-1}-n) V_2+
\frac{\vec\beta\cdot\dot V_1}{\alpha}(\vec\gamma-(n+A^{-1})\vec\beta)
\cr
+&\frac{\vec\gamma\cdot V_1}{\alpha}((n+A^{-1})\vec\beta-\vec\gamma)= n\left[\vec\xi[\nabla w]-\frac{p[\nabla w]-\vartheta}{\alpha}\vec\beta\right]
+\frac{p[\nabla w]-\vartheta}{\alpha}\vec\gamma
-\vec\zeta[\nabla w]-\vec f{\,}^\flat.
\end{aligned}
\end{equation}
Now (\ref{eq:u-energy}) implies that 
\[
\int_0^\infty \left(|V_1|^2+| \dot V_1|^2+|\nabla w|^2\right)\,e^{-nt}\,dt<\infty,
\]
and if we use the second equation in (\ref{eq:ODE6}) we see that
\[
\int_0^\infty \left(|V_2|^2+| \dot V_2|^2\right)\,e^{-nt}\,dt<\infty.
\]
We can summarize this as
\begin{equation}
\int_0^\infty \left(|V|^2+|\dot V|^2+|\nabla w|^2\right)\,e^{-nt}\,dt<\infty.
\label{eq:V-energy}
\end{equation}

Notice that the terms involving $\dot V_1$ and $\dot V_2$ in (\ref{eq:ODE6}) are of the form $(I+D(t))\dot V$ where $I$ is the identity matrix and  $\pmb{\bigr |}D(t)\pmb{\bigr |}\leq c\,\e^2(t)$. So we may multiply (\ref{eq:ODE6}) by $(I+D(t))^{-1}$ and,
after some calculations, see that $V$ satisfies a 1st-order system in the form
\begin{subequations}\label{ODE7:joint}
\begin{equation} \label{ODE7:a}
\frac{dV}{dt}+M(t)
V= F(t,\nabla w)+ F_0(t),
\end{equation}
where  $M(t)$ is a $2n\times 2n$ matrix of the form
\begin{equation}\label{ODE7:b}
\begin{aligned}
M(t)&=\, M_\infty +S_1(t)+S_2(t), \\
M_\infty&=\begin{pmatrix}
 -I & nI \\ (1-n^{-1})I & (1-n)I
 \end{pmatrix},\\
  S_1(t)&=
\begin{pmatrix}
I-A^{-1}B & A^{-1}-nI \\ C-BA^{-1}B+\frac{1-n}{n}I & BA^{-1}-I 
\end{pmatrix}, 
\end{aligned}
 \end{equation}
the $S_i$ satisfy
 \begin{equation}\label{ODE7:c}
 \begin{aligned}
\pmb{\bigr |}S_1(t)\pmb{\bigr |}\leq\e(t) \ \hbox{and}\  \pmb{\bigr |}S_2(t)\pmb{\bigr |}\leq c\,\e^2(t) \ \hbox{for}\ t>0,\\
 S_1(t)=0=S_2(t)  \ \hbox{for}\ t<0,
 \end{aligned}
 \end{equation}
the vector  $F(t,\nabla w)$ satisfies
 \begin{equation}\label{ODE7:d}
 | F(t,\nabla w)|\leq \ c\, \e(t)\,\meanint_{S^{n-1}} |\nabla w|\,ds\ \hbox{for}\ t>0
\ \hbox{and}\ 
F(t,\nabla w)\equiv 0\  \hbox{for}\ t<0,
 \end{equation}
 and $F_0(t)$ has support in $\log 2 \leq t \leq 2\log 2$ with $L^1$-norm satisfying
  \begin{equation}\label{ODE7:e}
 \| F_0\|_1\leq \ c\, (\|\vec f\|_p+\|f_0\|).
 \end{equation}
 \end{subequations}
 We can calculate the eigenvalues of $M_\infty$ to be $\lambda =0$ ($n$ times) and $\lambda =-n$ ($n$ times).
 The matrix
 \begin{subequations}\label{def-V:joint}
  \begin{equation}
J=
\begin{pmatrix} \label{def-V:a}
 nI & nI \\  I & (1-n)I
\end{pmatrix}
 \end{equation}
  diagonalizes $M_\infty$, i.e.\ 
\begin{equation} \label{def-V:b}
J^{-1}M_\infty J=\hbox{diag}(0,\dots,0,-n,\dots,-n),
 \end{equation}
so we introduce the change of dependent variables $V\to (\phi,\psi)$ by
 \begin{equation} \label{def-V:c}
V=J
\begin{pmatrix}
 \phi \\  \psi
\end{pmatrix}.
 \end{equation}
 \end{subequations}
 We find that the dynamical system (\ref{ODE7:a}) now takes the form
  (\ref{ourODEsystem:a}), where the conditions (\ref{ourODEsystem:c}) and  (\ref{ourODEsystem:d}) follow from 
  (\ref{ODE7:d}) and (\ref{ODE7:e}) respectively, 
  and ${\mathcal R}$ is of the form (\ref{DS:R}) 
  with
 \begin{equation}
R_1=\frac{n-1}{n^2}A^{-1}-\frac{n-1}{n}A^{-1}B+C-BA^{-1}B+\frac{1}{n}BA^{-1}-I.
\end{equation}
To simplify this expression for $R_1$, let us write
\[
A=n^{-1}(1+\widetilde A), 
\qquad
B=n^{-1}(1+\widetilde B), 
\qquad\hbox{and}\quad
C=n^{-1}(1+\widetilde C),
\]
where
$\pmb{|}\widetilde A\pmb{|}, \  \pmb{|}\widetilde B\pmb{|}, \  \pmb{|}\widetilde C\pmb{|} \leq c\,\e(t)$ as $t\to\infty$.
Then
\[
A^{-1}\approx n(I-\widetilde A),
\]
where $\approx$ means that the difference is bounded by $c\,\e^2(t)$ as $t\to\infty$, and a calculation shows
\[
R_1\approx \widetilde C-\widetilde B=C-nB\quad\hbox{as}\ t\to\infty.
\]
Since
\[
C-nB=\meanint_{S^{n-1}}(A-n\,A\,\theta\otimes\theta)\,ds_\theta,
\]
we see that (\ref{ourODEsystem:b}) holds and we have completed our derivation of (\ref{ourODEsystem:joint}).
\hfill$\Box$


\end{document}